\documentclass{amsart}
\usepackage{enumerate}

\usepackage[left=30mm, right=30mm, bottom=30mm]{geometry}
\usepackage{graphicx}
\usepackage{color}
\usepackage{mathtools}
\usepackage{tikz-cd}
\usepackage{amsmath,amssymb,amsfonts,amsthm}

\newtheorem{thm}{Theorem}[section]
\newtheorem*{thm*}{Theorem}
\newtheorem{prop}[thm]{Proposition}
\newtheorem{lemma}[thm]{Lemma}

\newtheorem{dfn}[thm]{Definition}

\theoremstyle{remark}
\newtheorem{remark}[thm]{Remark}

\newcommand{\C}{\mathbb{C}}
\newcommand{\F}{\mathbb{F}}
\newcommand{\Q}{\mathbb{Q}}
\newcommand{\R}{\mathbb{R}}

\newcommand{\Z}{\mathbb{Z}}

\newcommand{\Cbf}{\mathbf{C}}

\newcommand{\Disc}{\mathrm{Disc}\,}

\newcommand{\ep}{\varepsilon}

\newcommand{\con}{\equiv}

\newcommand{\ndiv}{\nmid}
\newcommand{\modd}[1]{\; ( \text{mod} \; #1)}

\newcommand{\maps}{\rightarrow}
\newcommand{\intersect}{\cap}

\newcommand{\union}{\cup}

\newcommand{\supp}{{\rm supp \;}}

\newcommand{\al}{\alpha}
\newcommand{\be}{\beta}

\newcommand{\ga}{\gamma}
\newcommand{\del}{\delta}

\newcommand{\om}{\omega}

\newcommand{\Ga}{\Gamma}

\newcommand{\Pcal}{\mathcal{P}}
\newcommand{\Bcal}{\mathcal{B}}

\newcommand{\Scal}{\mathcal{S}}

\newcommand{\Ubf}{\mathbf{U}}

\newcommand{\Xbf}{\mathbf{X}}

\newcommand{\zerobf}{\boldsymbol0}

\newcommand{\abf}{{\bf a}}

\newcommand{\kbf}{\mathbf{k}}

\newcommand{\ubf}{{\bf u}}
\newcommand{\vbf}{{\bf v}}

\newcommand{\x}{{\bf x}}
\newcommand{\xbf}{{\bf x}}

\newcommand{\beq}{\begin{equation}}
\newcommand{\eeq}{\end{equation}}

\makeatletter
\def\@tocline#1#2#3#4#5#6#7{\relax
  \ifnum #1>\c@tocdepth 
  \else
    \par \addpenalty\@secpenalty\addvspace{#2}%
    \begingroup \hyphenpenalty\@M
    \@ifempty{#4}{%
      \@tempdima\csname r@tocindent\number#1\endcsname\relax
    }{%
      \@tempdima#4\relax
    }%
    \parindent\z@ \leftskip#3\relax \advance\leftskip\@tempdima\relax
    \rightskip\@pnumwidth plus4em \parfillskip-\@pnumwidth
    #5\leavevmode\hskip-\@tempdima
      \ifcase #1
       \or\or \hskip 1em \or \hskip 2em \else \hskip 3em \fi%
      #6\nobreak\relax
    \hfill\hbox to\@pnumwidth{\@tocpagenum{#7}}\par
    \nobreak
    \endgroup
  \fi}
\makeatother

\theoremstyle{plain}

\DeclareMathOperator{\bfk}{{\mathbf{k}}}
\DeclareMathOperator{\bfa}{\mathbf{a}}
\DeclareMathOperator{\bfb}{\mathbf{b}}
\DeclareMathOperator{\bfe}{\mathbf{e}}
\DeclareMathOperator{\bfu}{\mathbf{u}}
\DeclareMathOperator{\bfv}{\mathbf{v}}

\DeclareMathOperator{\bfC}{\mathbf{C}}
\DeclareMathOperator{\bfU}{\mathbf{U}}
\DeclareMathOperator{\bfx}{\mathbf{x}}
\DeclareMathOperator{\X}{\mathbf{X}}
\DeclareMathOperator{\bX}{\mathbf{X}}
\DeclareMathOperator{\bx}{\mathbf{x}}
\DeclareMathOperator{\bu}{\mathbf{u}}
\DeclareMathOperator{\rank}{rank}
\DeclareMathOperator{\sing}{Sing}

\theoremstyle{definition}

\theoremstyle{remark}

\definecolor{pink}{rgb}{1,.2,.6}
\definecolor{orange}{rgb}{0.7,0.3,0}
\definecolor{blue}{rgb}{.2,.6,.75}
\definecolor{green}{rgb}{.4,.7,.4}
\definecolor{purple}{RGB}{127,0,255}

\numberwithin{equation}{section}

\begin{document}
 
\title[Application of a polynomial sieve]{Application of a polynomial sieve: \\beyond separation of variables}

\author{Dante Bonolis}
\address[Dante Bonolis]{Mathematics Department, Duke University,
120 Science Drive, Durham, North Carolina 27708, USA} \email{dante.bonolis@duke.edu}

\author{Lillian B. Pierce}
\address[Lillian B. Pierce]{
Mathematics Department,  Duke University,
120 Science Drive, Durham, North Carolina 27708, USA
} \email{pierce@math.duke.edu}

\begin{abstract}
Let a polynomial $f \in \Z[X_1,\ldots,X_n]$ be given. The square sieve can provide an upper bound for the number of integral  $\xbf \in [-B,B]^n$ such that $f(\xbf)$ is a perfect square. Recently this has been generalized substantially: first to a power sieve, counting $\xbf \in [-B,B]^n$ for which $f(\xbf)=y^r$ is solvable for $y \in \Z$; then to a polynomial sieve, counting $\xbf \in [-B,B]^n$ for which $f(\xbf)=g(y)$ is solvable, for a given polynomial $g$. 
Formally, a polynomial sieve lemma can encompass the more general problem of counting $\xbf \in [-B,B]^n$ for which $F(y,\xbf)=0$ is solvable, for a given polynomial $F$. Previous applications, however, have only succeeded in the case that $F(y,\xbf)$ exhibits separation of variables, that is, $F(y,\xbf)$ takes the form $f(\xbf) - g(y)$. In the present work, we present the first application of a polynomial sieve to count $\xbf \in [-B,B]^n$ such that $F(y,\xbf)=0$ is solvable, in a case for which $F$ does not exhibit separation of variables. Consequently, we obtain a new result toward a question of Serre, pertaining to counting points in thin sets. \\

\noindent
NOTE: Appended to the end of this paper, please find a correction, as published in the journal in which the original paper appeared. No changes have been made to the main body of the paper.
\end{abstract}
 
\maketitle

\section{Introduction}

Fix an integer $m \geq 2$ and integers $d, e \geq 1$. Consider the polynomial
\beq\label{F_dfn}
F(Y,\X)=Y^{md}+Y^{m(d-1)}f_{1}(\X)+\ldots+Y^mf_{d-1}(\X)+f_{d}(\X),
\eeq
in which for each $1 \leq i \leq d$,  $f_i \in \Z[X_1,\ldots, X_n]$ is a form with $\deg f_{i}=m\cdot e\cdot i$.
We are interested in counting
\[
N(F, B):=|\{\x\in [-B,B]^{n} \intersect \Z^n: \exists y\in\mathbb{Z}\text{ such that } F(y,\x)=0\}|.
\]
Trivially,   $N(F,B) \ll B^n$; our main result proves a nontrivial upper bound. 
We assume in what follows that $f_{d}\not\con 0$, since otherwise $(0,\Xbf)$ is a solution to $F(Y,\Xbf)=0$ for all $\Xbf$, and then $B^n \ll N(F,B)\ll B^n$. 
 (Throughout, we use the convention that $A \ll_\kappa B$ if there exists a constant $C$, possibly depending on $\kappa$, such that $|A| \leq C B.$)

\begin{thm}\label{thm_main}
Fix $n \geq 3$. Fix integers $m \geq 2$ and $e,d \geq 1$.
Let $F$ be defined as in (\ref{F_dfn}), with $f_d \not\con 0$. Suppose that  the weighted hypersurface $V(F(Y,\X)) \subset \mathbb{P}(e,1,\ldots,1)$ defined by $F(Y,\Xbf)=0$ is nonsingular over $\C$. Then 
\[ N(F,B) \ll B^{n-1 + \frac{1}{n+1}} (\log B)^{\frac{n}{n+1}}.
\]
  The implicit constant may depend on $n,m,d,e$, but is otherwise independent of $F$. 
\end{thm}
 
 The main progress achieved in Theorem \ref{thm_main} is for $n \geq 4$, $e \geq 2, d \geq 2$.
The requirement that $n \geq 3$ occurs since a key step, Proposition \ref{prop_G_no_linear}, is not true for $n=2$ (see Remark \ref{remark_badbad_n2}). In any case,   for $n=2,3$ the result of Theorem \ref{thm_main} is superceded by results of  Broberg in \cite{Bro03}, as described below in (\ref{Bro_n2}) and (\ref{Bro_n3}). 
When $e=1$, the variety $V(F(Y,\X)) \subset \mathbb{P}(e,1,\ldots,1)$ is unweighted, so that in the setting of Theorem \ref{thm_main}, to bound $N(F,B)$ it is equivalent to count points $[Y:X_1:\cdots:X_n]$ with $|Y|,|X_i| \ll B$ on a nonsingular projective hypersurface of degree at least 2 in $\mathbb{P}^n$. Then the result of Theorem \ref{thm_main} (in the stronger form   $N(F,B)\ll_{m,d,n,\ep} B^{n-1+\ep}$) has already been obtained by work of Heath-Brown and Browning, appearing in \cite{HB94,HB02,Bro03a,BroHB06a,BroHB06b}, as summarized by Salberger in  \cite{Sal07}.
Finally, when $d=1$, the result of Theorem \ref{thm_main} (aside from uniformity in the coefficients of $F$) follows from recent work of the first author in \cite{Bon21} (see Remark \ref{remark_d1}).
 
  The condition  $m\geq 2$ is applied in two ways: first, in the construction of certain sieve weights (see \S \ref{sec_lit} and the proof of Lemma \ref{lemma_sieve}), and second,   in \S \ref{sec_smooth} when we pass from the weighted variety to an unweighted variety. For illustration, we also describe how an alternative approach to the sieve lemma, conditional on GRH, can be devised when $m=1$ (see \S \ref{sec_GRH} and Remark \ref{remark_GRH}).

 Bounding $N(F,B)$ relates to a question of Serre on  counting integral points in thin sets. 
 Let $\mathcal{V}$ denote the affine variety 
\beq\label{Vcal}
\mathcal{V}= \{ (Y,\X) \in \mathbb{A}^{n+1} : F(Y,\X)=0\},
\eeq
and consider the projection 
\beq\label{proj_pi}
\pi:\begin{matrix}
\mathcal{V} &\rightarrow & \mathbb{A}^{n}\\
(y,\bfx) & \mapsto & \bfx.
\end{matrix}
\eeq
Under the hypotheses of Theorem \ref{thm_main},
  the set $Z=\pi(\mathcal{V}(\mathbb{Q}))$ is a \emph{thin set of type II} in $\mathbb{A}^n_\Q$, in the nomenclature of Serre. 
 Serre has posed a general question  that can be interpreted in our present setting as asking whether it is possible to prove that
\beq\label{Serre_question_0}
N(F,B) \ll B^{n-1}(\log B)^c
\eeq
for some $c$. Previous work by Broberg \cite{Bro03} nearly settled Serre's conjecture for thin sets of type II in $\mathbb{P}^{n-1}$  for $n=2,3$; see (\ref{Bro_n2}) and (\ref{Bro_n3}) below.
For $n \geq 4$, Theorem \ref{thm_main} represents new progress toward resolving Serre's question for certain thin sets of type II. Note that 
as $n \maps \infty$, the bound in Theorem \ref{thm_main} approaches a bound of the strength (\ref{Serre_question_0}).
We provide general background on Serre's question, and state precisely how Theorem \ref{thm_main} relates to previous literature on this question, in \S \ref{sec_Serre} and \S \ref{sec_lit}.

To prove Theorem \ref{thm_main}, we develop an appropriate polynomial sieve lemma, and then bound each contribution to the sieve using analytic, algebraic, and geometric ideas. A novel feature of this work is that we do not assume that $F(Y,\Xbf)$ exhibits separation of variables: that is, when $d \geq 2$, $F(Y,\Xbf)$ of the form (\ref{F_dfn})  cannot in general be written as $F(Y,\Xbf) = g(Y) - G(\Xbf)$ for polynomials $g,G$.  A formal polynomial sieve lemma has been formulated previously in a level of generality that does not require separation of variables;  see \cite{Bro15,BCLP23}. However,  in those works it has so far only been applied to count points on a variety that does exhibit separation of variables.
To our knowledge, Theorem \ref{thm_main} is the first application 
of a polynomial sieve to produce an upper bound for $N(F,B)$ in a case without separation of variables. We state precisely how Theorem \ref{thm_main} relates to previous literature on so-called square, power, and  polynomial sieves in \S \ref{sec_lit}.

A second strength of Theorem \ref{thm_main} is that the exponent in the upper bound for $N(F,B)$ is independent of $e$, where we recall that as a function of $\Xbf$, $F$ has highest degree $m \cdot e \cdot d$. For any given $\xbf \in [-B,B]^n$ such that $F(Y,\xbf)=0$ is solvable, one observes that any solution $y$ to $F(y,\xbf)=0$ must satisfy $y \ll B^e,$ and there can be at most $md$ solutions $y$ for the given $\xbf$ (or, equivalently, pre-images under the projection $\pi$ in (\ref{proj_pi})), since the coefficient of $Y^{md}$ in $F(Y,\Xbf)$ is nonzero. Thus an alternative method to bound $N(F,B)$ (up to an implicit constant depending on $md$) would be to count all $(n+1)$-tuples $\{(y,\xbf): y\ll B^e, x_i \ll B: F(y,\xbf)=0\}$.  Other potential methods might be sensitive to the role of $e$ or size of $d,m,$ (see for example Remark \ref{remark_e}), while in contrast both the method and the result of Theorem \ref{thm_main} do not depend on $e$ (aside from a possible implicit constant).

Third, we note that the result of Theorem \ref{thm_main} is independent of the coefficients of $F$; the implicit constant depends only on $F$ in terms of its degree. To accomplish this, we adapt a strategy of \cite{HB02}, also recently applied in a similar setting in \cite{BonBro22}, to show that either $N(F,B)$ is already acceptably small, or   $\|F\| \ll B^{(mde)^{n+2}}$. In the latter case, we then show that any dependence on $\|F\|$ in the sieve method is at most logarithmic, which we show is allowable for the result in Theorem \ref{thm_main}.

\subsection{Context of Theorem \ref{thm_main} within the study of Serre's question on thin sets}\label{sec_Serre}

Here we recall the notion of thin sets defined by Serre in \cite[\S 9.1 p.121]{Ser97} and \cite[p. 19]{Ser92}. 
Let $k$ be a field of characteristic zero and let $V$ be an irreducible  algebraic variety in $\mathbb{P}_k^n$ (respectively $\mathbb{A}_k^n$). A subset $M$ of $V(k)$ is said to be a projective (respectively, affine) thin set of type I if there is a closed subset $W \subset V$, $W \neq V$, with $M \subset W(k)$ (i.e. $M$ is not Zariski dense in $V$). A subset $M$ of $V(k)$ is said to be a projective (respectively, affine) thin set of type II if there is an irreducible projective (respectively, affine) algebraic variety $X$  with $\dim X = \dim V$, and a generically surjective morphism $\pi : X \maps V$  of degree $d \geq 2$ with $M \subset \pi(X(k))$.  Any thin set is a finite union of thin sets of type I and thin sets of type II. From now on we consider only $k=\Q$, although Serre's treatment considers any number field.

Given a thin set $M \subset \mathbb{A}^n_\Q$, define the counting function
\[ M(B) := | \{ \x \in M \intersect \Z^n : \max_{1 \leq i \leq n}|x_i| \leq B\}|,\]
so that trivially $M(B) \ll B^n$ for all $B \geq 1$. 
A theorem of Cohen \cite{Coh81} (see also \cite[Ch. 13 Thm. 1 p. 177]{Ser97}) shows that
\beq\label{Cohen_bound}
M(B) \ll_M B^{n-1/2} (\log B)^\ga \qquad \text{for some $\ga<1$,}\eeq
 where $\ll_M$ denotes that the implicit constant can depend on the coefficients of the equations defining $M$.
 As Serre remarks, this bound is essentially optimal, since the thin set 
 \beq\label{Cohen_set}
 M=\{\xbf =(x_1,\ldots,x_n)\in \Z^n : \text{$x_1$ is a square}\}\eeq
 has $M(B) \gg B^{n-1/2}.$ However, this $M$ arises from a morphism that is singular;
  it is reasonable to expect that the result can be improved under an appropriate nonsingularity assumption (such as in the setting of Theorem \ref{thm_main}).  
   
Now let  $M \subset \mathbb{P}_\Q^{n-1}$ be a thin set in projective space. Define the height function $H(x)$ for $x = [x_1:\ldots:x_n] \in  \mathbb{P}_\Q^{n-1}$ such that $(x_1,\ldots,x_n) \in \mathbb{Z}^n$ and $\gcd(x_1,\ldots,x_n)=1$ by $H(x) = \max_{1 \leq i \leq n} |x_i|$. Define the associated counting function 
\[ M_H(B) = \{ x \in M(\Q) : H(x) \leq B \}
\]
so that trivially $M_H(B) \ll B^{n}.$
Serre deduces in \cite[Ch. 13 Thm. 3]{Ser97} from an application of (\ref{Cohen_bound}) that 
\beq\label{Serre_Cohen_deduction}
M_H(B)
\ll_M B^{n-1/2} (\log B)^\ga \qquad \text{for some $\ga<1$.}
\eeq
   Serre raises a general question in \cite[p. 178]{Ser97}: is it possible to prove that 
\beq\label{Serre_question}
M_H(B) \ll B^{n-1}(\log B)^c
\eeq
for some $c$?
 (Note that the set (\ref{Cohen_set}) is not an example of a thin set here because if we set $M = \{[x_1^2:x_2: \cdots : x_n]\} \subset \mathbb{P}_\Q^{n-1}$ then for any $x_1 \neq 0$,
\[[x_1:x_2: \cdots : x_n] = x_1[x_1:x_2: \cdots : x_n]  = [x_1^2:x_1x_2: \cdots : x_1x_n] \in M 
\]
so that $M \supset \mathbb{P}_\Q^{n-1}$.)

\subsubsection{Results for thin sets of type I} 
If $Z$ is an  irreducible projective variety in $\mathbb{P}_\Q^{n-1}$ of degree $d \geq 2$, Serre deduces from (\ref{Serre_Cohen_deduction}) that $Z_H(B) \ll_Z B^{\dim Z+1/2} (\log B)^\ga$ for some $\ga<1$. Serre   asks if it is possible to prove that $Z_H(B) \ll_Z B^{\dim Z} (\log B)^c$ for some $c$.  
(This question is raised in both \cite[p. 178]{Ser97} and \cite[p. 27]{Ser92}. Serre provides an example of a quadric for  which a logarithmic factor necessarily arises. See also the question in the case of a hypersurface in Heath-Brown \cite[p. 227]{HB83}, formally stated in both non-uniform and uniform versions as \cite[Conj. 1, Conj. 2]{HB02}.)    This is now called the \emph{dimension growth conjecture} (in the terminology of \cite{Bro09}), and is often described as the statement that  \beq\label{DimGrowth}
Z_H(B)\ll_{Z,\ep} B^{\dim Z+\ep} \qquad \text{for every $\ep>0$.}
\eeq
A refined version, credited to Heath-Brown and known as the \emph{uniform dimension growth conjecture}, is the statement that \beq\label{UDimGrowth}
Z_H(B)\ll_{n,\deg Z,\ep} B^{\dim Z+\ep} \qquad \text{for every $\ep>0$.}
\eeq

 In the case that $Z \subset \mathbb{P}_\Q^{n-1}$ is a nonsingular projective hypersurface of degree $d \geq 2$, as mentioned before, combined works of Browning and Heath-Brown have proved (\ref{UDimGrowth}) for all $n \geq 3.$ More generally, Browning, Heath-Brown and Salberger proved (\ref{UDimGrowth}) for all   geometrically integral varieties of degree $d=2$ and $d \geq 6$ (see \cite{HB02} and \cite{BHS06}, respectively). Recent work of Salberger has proved (\ref{DimGrowth}) in all remaining cases, and has even proved the uniform version (\ref{UDimGrowth}) for $d \geq 4$ \cite{Sal23}. 
 See \cite{CCDN20} for a helpful survey, statements of open questions, and  new progress such as an explicit bound $Z_H(B) \leq C d^{E} B^{\dim Z}$ when $\deg Z =d\geq 5$, for a certain $C=C(n)$ and $E=E(n).$
   The resolution of the dimension growth conjecture means that attention now turns to thin sets of type II, the subject of the present article.

\subsubsection{Results for thin sets of type II} We turn  to the case of thin sets of type II, our present focus. 
 Given a finite cover $\phi : X \maps \mathbb{P}^{n-1}$ over $\Q$ with $n \geq 2$, $X$ irreducible  and $\phi$ of degree at least 2, set 
\beq\label{NBphi}
N_B(\phi) = |\{ P \in X(\Q): H(\phi(P)) \leq B\}|
\eeq
for the standard height function above. Serre's question asks whether \beq\label{Serre_again}
N_B(\phi)\ll_{\phi,n} B^{n-1}(\log B)^c \qquad \text{ for some $c$,}
\eeq
or in a uniform version, 
\beq\label{Serre_again_unif}
N_B(\phi)\ll_{\deg \phi, n} B^{n-1}(\log B)^c \qquad \text{ for some $c$.}
\eeq

For $n=2,3$ work of Broberg  via the determinant method proves   cases of Serre's conjecture  up to the logarithmic factor \cite{Bro03}. Precisely, for $\phi: X \maps \mathbb{P}^1$ of degree $r \geq 2$, Broberg proves 
 \beq\label{Bro_n2}
 N_B(\phi)\ll_{\phi,\ep} B^{2/r+\ep} \qquad \text{for any $\ep>0$.}
 \eeq
 For $\phi: X \maps \mathbb{P}^2$ of degree $r$, Broberg proves 
 \beq\label{Bro_n3}
 \text{$N_B(\phi)\ll_{\phi,\ep} B^{2+\ep}$ for $r \geq 3$,   \qquad $N_B(\phi)\ll_{\phi,\ep} B^{9/4+\ep}$ for $r=2$, for any $\ep>0$.}
 \eeq
 For $n \geq 4,$ the question remains open whether one can achieve $N_B(\phi)\ll B^{n-1+\ep}$ for all $\ep>0$, although we record some progress on this for specific types of $\phi$ in \S \ref{sec_lit}. 
 
Now recall the setting of Theorem \ref{thm_main} in this paper, and the affine variety $\mathcal{V} \subset \mathbb{A}^{n+1}$ defined in (\ref{Vcal}) according to the polynomial $F(Y,\X)$. Under the hypotheses of Theorem \ref{thm_main}, we have:
\begin{itemize}
    \item[$i)$] The variety $\mathcal{V}$ is irreducible (see Remark \ref{remark_V_irred});
    \item[$ii)$] The projection $\pi$ has degree $dm>1$ since $m \geq 2$.
\end{itemize}

Thus  $Z=\pi(\mathcal{V}(\mathbb{Q}))$ is a thin set of type II in $\mathbb{A}^n_\Q$, and in particular Cohen's result (\ref{Cohen_bound}) implies that 
\beq\label{Cohen} Z(B)=N(F,B) \ll_F B^{n-1/2} (\log B)^\ga,
\eeq
following the same reasoning as \cite[Ch. 13 Thm. 2 p. 178]{Ser97}. 
Or, interpreting the setting of Theorem \ref{thm_main} as counting points on a finite cover $\phi$ of $\mathbb{P}^{n-1}$ as in (\ref{NBphi}), this shows
\[ N_B(\phi) \ll N(F,B) \ll_\phi B^{n-1/2}( \log B)^\ga.\]
 Our new work,  Theorem \ref{thm_main}, improves on (\ref{Cohen}) for each $n \geq 3$, for $F$ of the form (\ref{F_dfn}) with $V(F(Y,\X))$ nonsingular, and approaches a uniform bound of the strength (\ref{Serre_again_unif}) as $n \maps \infty$.

\subsection{Context of Theorem \ref{thm_main} within sieve methods}\label{sec_lit}
We now recall a few recent developments of sieve methods in the context of counting solutions to Diophantine equations, with a particular focus on progress toward Serre's conjecture for type II sets, as described above. 

\subsubsection{Square sieve}
Let $f(\Xbf) \in \Z[X_1,\ldots,X_n]$ be a fixed polynomial. Let  $\Bcal$ be a ``box,'' such as $ [-B,B]^n$ or more generally $\prod_i [-B_i, B_i].$ In \cite{HB84}, Heath-Brown codified the square sieve to count the number of integral values $\xbf \in \Bcal$ such that $f(\xbf)=y^2$ is solvable over $\Z$, building on a method of Hooley \cite{Hoo78}. At its heart was a formal sieve lemma involving a character sum with Legendre symbols.  Heath-Brown applied this in particular to improve the error term in an asymptotic for the number of consecutive square-free numbers in a range. 
In \cite{Pie06}, Pierce developed a stronger version of the square sieve, with a sieving set comprised of products of two primes rather than primes; this effectively allows the underlying modulus to be larger relative to the box $\Bcal$, by factoring the modulus and using the $q$-analogue of van der Corput differencing. Pierce applied this to prove a nontrivial upper bound for $3$-torsion in class groups of quadratic fields \cite{Pie06};  Heath-Brown subsequently used this sieve method to prove there are finitely many imaginary quadratic fields having class group of exponent 5 \cite{HB08}; Bonolis and Browning applied it to prove a uniform bound for counting rational points on hyperelliptic fibrations
\cite{BonBro22}.

\subsubsection{Power sieve}\label{sec_power_Serre}

 The square sieve has been generalized to a power sieve, in order to count integral values $\xbf \in \Bcal$ with $f(\xbf)=y^r$ solvable, for a fixed $r \geq 2$. Recall the question of bounding $N_B(\phi)$  as in (\ref{Serre_again}).   For any $n \geq 2$, in the special case that $\phi$ is a nonsingular cyclic cover of degree $r \geq 2,$ Munshi observed this can be reduced to counting the number of integral values   $\xbf \in [-B,B]^n$ with $F(x_1,\ldots,x_n)=y^r$ solvable, for a nonsingular form $F$ of degree $mr$ for some $m\geq 1$. 
To bound this, Munshi developed  a formal sieve lemma involving a character sum in terms of multiplicative Dirichlet characters \cite{Mun09}.
Munshi applied it to prove that 
\beq\label{Munshi}
|\{\xbf \in [-B,B]^n : \text{$F(\xbf)=y^r$ is solvable over $\Z$}\}| \ll B^{n-1+ \frac{1}{n}}(\log B)^{\frac{n-1}{n}}
\eeq
Consequently, this proved $N_B(\phi)\ll B^{n-1+ \frac{1}{n}} (\log B)^{\frac{n-1}{n}}$
for nonsingular cyclic covers.
(See \cite[Remark 1]{Bon21} for a note on the history of this result; the exponents stated here are slightly different from those presented in \cite{Mun09}.)  

In \cite{HBPie12} Heath-Brown and Pierce have strengthened the power sieve, by using a sieving set comprised of products of primes, generalizing the approach of \cite{Pie06}. They used this method to prove that for any polynomial $f(\Xbf) \in \Z[X_1,\ldots,X_n] $ of degree $d \geq 3$ with nonsingular leading form, and for any $r \geq 2$, 
\beq\label{HBP_eqn}
|\{\xbf \in [-B,B]^n : \text{$f(\xbf)=y^r$ is solvable over $\Z$}\}|\ll 
\begin{cases}
    B^{n-1+\frac{n(8-n)+4}{6n+4}}(\log B)^2, & 2 \leq n \leq 8 \\
    B^{n-1+ \frac{1}{2n+10}}(\log B)^2, & n=9\\
    B^{n-1 - \frac{(n-10)}{2n+10}}(\log B)^2, & n \geq 10.
\end{cases}
\eeq
This  proves Serre's conjecture (\ref{Serre_again}) for $N_B(\phi)$, for all nonsingular cyclic covers, for $n \geq 10.$  
Indeed, the bound achieved is even smaller than the general conjecture, which is reasonable due to the imposed nonsingularity assumption.

 Independently, Brandes also developed a power sieve in \cite{Bra15}, applied to counting sums and differences of power-free numbers.
 
\subsubsection{Polynomial sieve: with separation of variables}
The next significant generalization addressed counting $\xbf \in \Bcal$ for which $g(y) = f(\xbf)$ is solvable, for appropriate polynomials $g,f$. Here, a quite general framework for a polynomial sieve lemma was developed by Browning  in \cite{Bro15}. Specifically, in that work, Browning applied the polynomial sieve lemma to count $x_1,x_2$ such that $g(y) = f(x_1,x_2)$ is solvable, for particular functions $f, g$, that enabled an application showing the sparsity of like sums of a quartic polynomial of one variable.

 Bonolis \cite{Bon21} further developed a polynomial sieve lemma with a character sum involving trace functions. Applying this, he proved that for any polynomial $g \in \Z[Y]$ of degree $r \geq 2,$ and any irreducible form $F \in \Z[X_1,\ldots, X_n]$ of degree $e \geq 2$ such that the projective hypersurface $V(F)$ defined by $F=0$ is nonsingular over $\C$, then 
 \beq\label{Bonolis}
|\{\xbf \in [-B,B]^n : \text{$F(\xbf)=g(y)$ is solvable over $\Z$}\}|\ll B^{n-1 + \frac{1}{n+1}}(\log B)^{\frac{n}{n+1}}.
\eeq
(This improves (\ref{Munshi}) and recovers the result initially stated in \cite{Mun09}; see \cite[Remark 1]{Bon21}.) 
 This can also be seen as an improvement on    Cohen's theorem (\ref{Cohen}) for a special type of thin set (defined as the image of $\mathcal{V}=\{(y,\xbf) \in \mathbb{A}^{n+1}: F(\xbf) - g(y)=0\}$ under $(y,\xbf) \mapsto \xbf$, under the assumption that $V(F)$ defines a nonsingular projective hypersurface).
 The special case of our Theorem \ref{thm_main} when $d=1$ follows from \cite[Theorem 1.1]{Bon21}; see Remark \ref{remark_d1}.

Notably, the method employed in \cite{Bon21} to prove (\ref{Bonolis}) was the first to demonstrate  nontrivial averaging over pairs of primes in the sieving set,  and exploiting such a strategy is  central to the strength of our main theorem. We explain explicitly the advantage of such averaging  in equations  (\ref{T_contribution}) and (\ref{upq}), below. For now, we simply state abstractly that any polynomial sieve method tests the solvability of the desired equation modulo $p$ for primes in a chosen sieving set $\mathcal{P}$. The outcome of applying a sieve lemma (such as Lemma \ref{lemma_sieve} below) is that one must bound from above an expression roughly of the form  $|\mathcal{P}|^{-2}\sum_{p \neq q \in \mathcal{P}} T(p,q)$, where $T(p,q)$  studies the solvability of the desired equation modulo pairs $p \neq q \in \mathcal{P}$. Previous to \cite{Bon21}, papers applying any type of polynomial sieve produced an upper bound for $|T(p,q)|$ that was uniform over $p,q$ and then summed trivially over $ p \neq q \in \mathcal{P}$. Instead, averaging nontrivially over $p,q$ exploits the fact that $T(p,q)$ is typically smaller than its worst (largest) upper bound.

 Most recently, a geometric generalization of Browning's polynomial sieve lemma has been developed  over function fields by Bucur, Cojocaru, Lal\'in and the second author in \cite{BCLP23}. They pose an analogue of Serre's question (\ref{Serre_question}) in that setting (also raised by Browning and Vishe \cite{BroVis15}), and apply   a polynomial sieve to prove a bound of analogous strength to  (\ref{Bonolis}), in the special case of nonsingular cyclic covers in a function field setting.  It remains an interesting open question to achieve a stronger bound such as (\ref{HBP_eqn}), or to prove results for finite covers that are noncyclic, in such a function field setting.

\subsubsection{Polynomial sieve: without separation of variables}\label{sec_intro_poly_sep}
So far we have mentioned applications of a sieve lemma to count solutions to $G(Y,\Xbf)=0$ when $G$ separates variables as $G(Y,\Xbf)=g(Y) - f(\Xbf)$ for some polynomials $g,f.$
 More generally, it is reasonable to ask---and this is a motivation for the present paper---whether an appropriate polynomial sieve can be employed to count solutions to equations of the form 
$G(Y,\Xbf)=0$ where $G(Y,\mathbf{X})\in\mathbb{Z}[Y,X_1,\ldots,X_n]$ is a polynomial of degree $D$ of the form
\beq\label{G_dfn_intro}
G(Y,\X)=Y^{D}+Y^{D-1}f_{1}(\X)+\ldots + Yf_{D-1}(\X)+ f_{D}(\X),
\eeq
where each $f_i$ is a form of degree $i \cdot e$, and we assume that the weighted hypersurface
$
V(G(Y,\X)))\subset\mathbb{P}(e,1,\ldots,1 )
$
defined by $G(Y,\X)=0$
is nonsingular.
Define  
\[
N(G, B):=|\{\x\in [-B,B]^{n}: \exists y\in\mathbb{Z}\text{ such that } G(y,\x)=0\}|.
\]
 Under the assumption $f_D \not\con 0$, the aim is to improve on the trivial bound $N(G,B) \ll B^n$.
To be clear, the formal sieve lemmas appearing in \cite{Bro15,BCLP23} include this level of generality, but have only been applied to prove a bound for $N(G,B)$  when separation of variables occurs.
  In this paper we accomplish  the first application of the polynomial sieve without assuming separation of variables, but under the additional 
assumption that the degree $D$ of $G(Y,\X)$ defined in (\ref{G_dfn_intro}) factors as $D=md$ for some $m \geq 2$, and all powers of $Y$ that appear are divisible by $m$.  (To see why this restriction is useful, see the proof of Lemma \ref{lemma_sieve}; for an alternative approach when $m=1$, conditional on GRH, see Remark \ref{remark_GRH} and \S \ref{sec_GRH}.)

 The strength of our approach hinges on a  particular formulation of the polynomial sieve, given in Lemma \ref{lemma_sieve}. It is worthwhile to compare our formulation   with the polynomial sieve presented in \cite[Theorem 1.1]{Bro15}. 
In \cite[Theorem 1.1]{Bro15}, the sieve weight system, adapted to counting solutions to (\ref{G_dfn_intro}), is defined as follows:
\[
w_{p,\text{Bro}}(\bfk) = \alpha + (\nu_{p}(\bfk) - 1)(D - \nu_{p}(\bfk))
,\]
in which $\nu_p (\kbf)= |\{ y \in \F_p: G(y,\kbf)=0 \in \F_p\}|.$
(These weights are then applied in an inequality analogous to (\ref{sieve_lemma_square}) below, to derive a sieve lemma.)
Consequently, if $G(Y,\bfk)=0$ is solvable over $\mathbb{Z}$, the conditions $1 \leq \nu_{p}(\bfk) \leq D$ and $\alpha > 0$ guarantee that $w_{p,\text{Bro}}(\mathbf{k}) > 0$ \textit{for any} $p$. In our approach, we consider simpler weights:
\[
w_{p}(\bfk) = \nu_{p}(\bfk) - 1.
\]
Thus, in our situation, if $G(Y,\bfk)=0$ is solvable over $\mathbb{Z}$, we can only conclude that $w_{p}(\bfk) \geq 0$. However, it is still possible to establish that $w_{p}(\mathbf{k}) > 0$ for \textit{a positive proportion of primes}, which suffices for our application. (Precisely, we obtain $\om_p(\kbf)>0$ for those $p \con 1 \modd{m}$ where $m \geq 2$; see (\ref{Pm}) in the proof of Lemma \ref{lemma_sieve}.)

The simplicity of our weight system turns out to be crucial for bounding the terms that appear in the polynomial sieve lemma. In the setting of the polynomial $F(Y,\X)$ as in (\ref{F_dfn}), our main task will be to prove square root cancellation for the sum
\[
\sum_{\substack{(z,\bfa ) \in\mathbb{F}_{p}^{n+1}\\ F(z^{e},\bfa)=0}}e_{p}(\langle\bfa,\bfu\rangle),
\]
for generic $\bfa\in\mathbb{F}_{p}^{n}$, which can be accomplished by exploiting the smoothness of the variety $V(F(Z^e,\X))$. On the other hand, if we were to adopt \cite[Theorem 1.1]{Bro15}, the presence of the factor $(\nu_{p}(\bfk))^{2}$ would lead to the exponential sum
\[
\sum_{\substack{(z_{1},z_{2},\bfa ) \in\mathbb{F}_{p}^{n+2}\\ F(z_{1}^{e},\bfa)=0\\F(z_{2}^{e},\bfa)=0}}e_{p}(\langle\bfa,\bfu\rangle),
\]
which is more challenging to handle, due to the highly singular nature of the variety $V(F(Z_{1}^e,\X)) \cap V(F(Z_{2}^e,\X)).$

\subsection{Overview of the method}\label{sec_overview}
   We now provide an overview of our method,   highlighting four key aspects of our  strategy.
  To prove a nontrivial upper bound for $N(F,B)$ via a sieve, we introduce a smooth non-negative function $W:\mathbb{R}^{n}\rightarrow\mathbb{R}_{\geq 0}$ defined by $W(\bfx) = w(\bfx/B)$, where $w$ is an infinitely differentiable, compactly supported function that is $\equiv 1$ on $[-1,1]^n$, and supported in $[-2,2]^{n}$. Define the smoothed counting function
\beq\label{S_count}
\mathcal{S}(F,B):=\sum_{\substack{\bfk\in\mathbb{Z}^{n}\\ F(y,\bfk)=0 \text{ solvable}}}W(\bfk),
\eeq
which sums over $\bfk\in \mathbb{Z}^n$ such that there exists $y \in \mathbb{Z}$ with $F(y,\bfk)=0$. By construction 
\[ N(F,B) \leq \mathcal{S}(F,B),\]
and we may focus on proving a nontrivial upper bound for $\mathcal{S}(F,B).$ 
We employ the following sieve lemma, which we prove in \S \ref{sec_sieve_proof}. Here and throughout, given a polynomial $f$, we let $\|f\|$ denote the maximum absolute value of any coefficient of $f$.
\begin{lemma}[Polynomial sieve lemma]\label{lemma_sieve}
Let $e,d\geq 1$ and  $m \geq 2$ be integers.
Consider the polynomial
\[
F(Y,\X)=Y^{md}+Y^{m(d-1)}f_{1}(\X)+\ldots+Y^m f_{d-1}(\X)+f_{d}(\X),
\]
under the assumption that $f_{d}\not\con 0$, and that $\deg f_{i}=m\cdot e\cdot i$ for each $1 \leq i \leq d$. 

Let $B \geq 1$ and define a smooth weight $W$ supported in $[-2B,2B]^n$ and $\con 1$ on $[-B,B]^n,$ as above. Let $\mathcal{P}\subset\{p\equiv 1\mod m\}$ be a finite set of primes $p \in [Q,2Q]$, with   cardinality $P$.  Suppose that $Q = B^{\kappa}$ for some fixed $0< \kappa \leq 1$ and that $P \gg Q/\log Q$. Suppose also that  
\beq\label{P_lower_lemma} P \gg_{m,e,d} \max\{ \log \|f_d\|, \log B\}.
\eeq
For each $\kbf \in \Z^n$ and $p \in \Pcal$ define 
\[\nu_p(\bfk) = |\{  y \in \F_p : F(y,\bfk)=0 \modd{p}\}|.\]
Then
\[
\mathcal{S}(F,B)\ll_{m,e,d} \sum_{\bfk:f_{d}(\bfk)=0}W(\bfk)+ \frac{1}{P}\sum_{\bfk}W(\bfk)+\frac{1}{P^{2}}\sum_{\substack{p,q\in\mathcal{P}\\p\neq q}}\left|\sum_{\bfk}W(\bfk)(\nu_{p}(\bfk)-1)(\nu_{q}(\bfk)-1)\right|.
\]
\end{lemma}
 \begin{remark}\label{remark_GRH}
We observe that the same lemma holds for  $m =1,$ conditional on GRH, with (\ref{P_lower_lemma}) replaced by 
$Q \gg_{m,e,d} \max \{(\log \|F\|)^{\al_0}, (\log B)^{\al_0}\}$ for some $\al_0>2.$
For the sake of illustration, we demonstrate this in \S \ref{sec_GRH}, although we do not apply such a conditional result in this paper.
\end{remark}
 
We now point out four key aspects of our method for applying this sieve lemma to prove Theorem \ref{thm_main}. First,   for all $\kbf$ and for all primes $p$, $\nu_p(\kbf) \leq md$; this is because $Y^{md}$ has coefficient 1 in $F(Y,\Xbf)$, so that  for all values of $\bfk$, $F(Y,\kbf)$ is of degree $md$ as a polynomial in $Y$. On the other hand, in the proof of the lemma, we use the assumption that each prime in the sieving set has $p \con 1 \modd{m}$ in order to provide a lower bound $\nu_p(\kbf)-1 \geq m-1>0$ for many $\kbf$, motivating our requirement that $m \geq 2.$  This is the first novelty of our method for dealing with a case in which the variables $Y,\Xbf$ are not ``separated.''

For each pair of primes $p \neq q \in \Pcal,$ the sieve lemma leads us to study
\beq\label{T_dfn}
T(p,q):=\sum_{\bfk \in \Z^n }W(\bfk)(\nu_{p}(\bfk)-1)(\nu_{q}(\bfk)-1).
\eeq
 After an application of the Poisson summation formula,  we see that
\[
T(p,q)=\left(\frac{1}{pq}\right)^{n}\sum_{\bfu \in \Z^n}\hat{W}\left(\frac{\bfu}{pq}\right) g(\bfu,pq),
\]
where
\begin{equation}
g(\bfu,pq):=\sum_{\bfa \modd{pq}}(\nu_{p}(\bfa)-1)(\nu_{q}(\bfa)-1)e_{pq}(\langle\bfa,\bfu\rangle).
\label{eq : expsumg}
\end{equation}
Here we write each coordinate of $\bfa$ in terms of its residue class modulo $pq$, and $e_{pq}(t) = e^{2\pi i t/pq}.$
After showing that $g(\ubf,pq)$ satisfies a multiplicativity relation, we can focus on the case of prime modulus, and study 
\[
g(\bfu, p):=\sum_{\bfa \in\mathbb{F}_{p}^{n}}(\nu_{p}(\bfa)-1)e_{p}(\langle\bfa,\bfu\rangle).
\]
We show that the main task to bound $g(\ubf,p)$ is to bound the exponential sum
\[
\sum_{\substack{(y,\bfa ) \in\mathbb{F}_{p}^{n+1}\\ F(y,\bfa)=0}}e_{p}(\langle\bfa,\bfu\rangle).
\]
Here we highlight a second aspect: the fact that the polynomial $F(Y,\X)$ is not homogeneous  motivates a more sophisticated approach  to bounding this sum (see Remark \ref{remark_F_hom}). Given a polynomial $H$, let $V(H)$ denote the corresponding variety $\{ H=0\}$, and let $\langle \Xbf,\Ubf \rangle = \sum_i X_i U_i.$ Roughly speaking, for each prime $p$ we divide  $\ubf \in \Z^n$ into three cases: a \emph{type zero} case when $\ubf \con 0 \modd{p}$,  a \emph{good} case when $V(\langle \Xbf,\ubf\rangle)$ is not tangent to $V(F(Y,\Xbf))$ over $\overline{\F}_p$, and finally a \emph{bad} case in which $V(\langle \Xbf,\ubf\rangle)$ is  tangent to $V(F(Y,\Xbf))$ over $\overline{\F}_p$.  
(More precisely, we reformulate this in terms of varieties in unweighted projective space.)
In the type zero case, we can only show that $g(\zerobf,p)\ll p^{n-1/2}$, but such cases are sparse. In the remaining two cases, we apply a version of the Weil bound to   $g(\bfu,p)$, obtaining $g(\ubf,p)\ll p^{n/2}$ if $\ubf$ is good and $g(\ubf,p)\ll p^{n/2+1/2}$ if $\ubf$ is bad (Proposition \ref{prop_g_sum}).

A third crucial aspect arises when we assemble this information efficiently inside the third term on the right-hand side of the sieve lemma, namely
\beq\label{T_contribution}
\frac{1}{P^2}\sum_{p \neq q \in \Pcal}|T(p,q)| \ll \frac{1}{P^2Q^{2n}}\sum_{p \neq q \in \Pcal} \sum_{\ubf \in \Z^n}\left| \hat{W}\left( \frac{\ubf}{pq}\right) g(\ubf,pq)\right|.
\eeq
 In many earlier applications of the power sieve or polynomial sieve to count solutions to Diophantine equations, the strategy has been to bound $|T(p,q)|$ uniformly over $p \neq q$ and simply sum trivially over $p \neq q$. However, recent work of the first author demonstrated how to take advantage of nontrivial averaging over the sum of $p \neq q \in \Pcal$; see \cite{Bon21}. In this paper, we also average nontrivially over $p \neq q$ and this contributes to  the strength of our main theorem.

In order to average nontrivially over $p \neq q \in \Pcal$, we quantify the fact that there cannot be many triples $\ubf,p,q$ for which $\ubf$ is simultaneously bad for both $p$ and  $q$. Roughly speaking, we characterize the dual variety of the original hypersurface $V(F(Y,\Xbf))$ according to an irreducible polynomial $G(U_Y,U_1,\ldots,U_n)$, and observe that $G(0,\ubf) \neq 0$ precisely when the hyperplane $V(\langle \ubf,\Xbf \rangle)$ is not tangent to $V(F(Y,\Xbf))$ over $\C$. Then we reverse the order of summation in the right-hand side of (\ref{T_contribution}), writing it as
\beq\label{upq}
\frac{1}{P^2Q^{2n}}\sum_{\ubf \in \Z^n} \sum_{p \neq q \in \Pcal} \left| \hat{W}\left( \frac{\ubf}{pq}\right) g(\ubf,pq)\right|.
\eeq
The sum over $\ubf$ can be split into case (a) where $G(0,\ubf) \neq 0$ and case (b) where $G(0,\ubf) = 0.$ In case (a), we show $\ubf$ is bad modulo $p$ and $q$ only if $p$ and $q$ divide the (nonzero) value of a certain resultant polynomial; thus there can only be very few such $p,q$. 

A fourth key aspect arises in case (b), for which $\ubf$ is bad for all primes (since the value of the resultant is zero). To compensate, we show that there are not too many $\ubf$ for which   $G(0,\ubf)=0.$ This step  is one of the significant novelties of the paper. It requires understanding not the  variety $V(G(U_Y, \Ubf))$ but $V(G(U_Y,\Ubf)) \intersect V(U_Y)$, the intersection  with the hyperplane $U_Y=0.$ To tackle this,  we show that any polynomial divisor of $G(0,\Ubf)$ has degree at least 2 (Proposition \ref{prop_G_no_linear}), so that we can apply strong bounds of Heath-Brown \cite{HB02} and Pila \cite{Pil95} to count solutions to $G(0,\ubf)=0$ (see (\ref{HBP_bounds})). To prove the key result in Proposition \ref{prop_G_no_linear},  we employ a geometric argument to show that given a nonsingular projective hypersurface $X$ and a projective line $\ell$ not contained in $X$, the generic hyperplane containing $\ell$ is not tangent to $X$. This statement, proved in \S \ref{sec_G} via a strategy suggested by Per Salberger, is critical to the method and the ultimate strength of Theorem \ref{thm_main}.

\begin{remark}\label{remark_e}
 
 It would be interesting to consider bounding $N(F,B)$, in the setting of Theorem \ref{thm_main}, by other methods. As mentioned earlier, one approach is to count all $(n+1)$-tuples $\{ (y,\x) \in \Z^{n+1} : y \ll B^e, x_i \ll B : F(y,\xbf)=0\},$  for example, by  applying the determinant method. Since the range of $y$ depends on $e$, such a direct approach is likely to produce a bound for $N(F,B)$ with an exponent depending on $e$. Alternatively,  
one could fix $x_2,\ldots,x_n$ (with $\approx B^{n-1}$ such choices) and consider the resulting equation as a projective curve in variables $y, x_1$. Supposing that the resulting curve    is generically of degree $dme$, an application of Bombieri-Pila \cite{BomPil89} could count $(y,x_1)$ in the square $[-B^e, B^e]^2$. This could ultimately lead to a total bound of the form $N(F,B)\ll B^{n-1} \cdot B^{e/dme+\ep} = B^{n-1+1/dm+\ep}$. This putative outcome appears independent of $e$, but the method has overcounted $x_1$ in the range $B^e$; nevertheless, such an approach could be advantageous for large $d,m.$
\end{remark}

\subsection{Notation}
 We use $e_q(t) = e^{2\pi i t/q}.$
 We denote   $\X=(X_{1},\ldots,X_{n})$, $\bfU=(U_{1},...,U_{n})$. Moreover, for two vectors $\mathbf{s}=(s_{1},\ldots s_{n}),\mathbf{t}=(t_{1},\ldots,t_{n})$, we define $\langle\mathbf{s},\mathbf{t}\rangle=\sum_{i=1}^{n}s_{i}t_{i}$.
 We let $\|F\|$ denote the absolute value of the maximum coefficient in a polynomial $F \in \Z[X_1,\ldots,X_n]$; similarly $\|\Xbf\| = \max_{1 \leq i \leq n}|X_i|$ for $\Xbf \in \Z^n$.

\subsection*{Acknowledgements}
 
The authors thank T. Browning for suggesting the application of the polynomial sieve to smooth coverings and for useful discussions, and J. Lyczak for many helpful remarks. In addition, the authors thank P. Salberger for suggesting a strategy to prove Proposition \ref{prop_G_no_linear}, and both Salberger and an anonymous referee for helpful remarks on an earlier version of this manuscript.
 The authors credit ChatGPT with expository edits to the last two paragraphs of \S \ref{sec_intro_poly_sep}.

\section{Reduction to remove dependence on $\|F\|$}
Recall that Theorem \ref{thm_main} states that the upper bound for $N(F,B)$ is only dependent on the degree of $F$, and not on the coefficients of $F$. In fact, the sieve methods we apply  prove an upper bound for $N(F,B)$ that \emph{can} depend on $\|F\|$. 
In this section we show by alternative methods that we may assume that $\|F\|\ll B^{(mde)^{n+2}}$. The method does not rely on assuming $m \geq 2$ in (\ref{F_dfn}), and so without any additional trouble we may work more generally in the setting of (\ref{G_dfn_intro}).
\begin{lemma}\label{lemma_F_small}
Let $V(G(Y,\X))\subset\mathbb{P}(e,1,\ldots,1)$ be defined by  
\[G(Y,\Xbf)=Y^D + Y^{D-1}f_1(\Xbf) + \cdots + Yf_{D-1}(\X)+ f_D(\Xbf)\]
with each $f_i$ a form of $\deg f_i = i \cdot e$, for fixed $D,e \geq 1$ and $n \geq 1$.  Assume that $f_D \not\con 0$ and the weighted hypersurface $V(G(Y,\X))\subset\mathbb{P}(e,1,\ldots,1)$ is absolutely irreducible. Then either
\[
\|G\|\ll B^{(De)^{n+2}},
\]
or $N(G,B)\ll_{n,D,e} B^{n-1}$. 
\end{lemma}
\begin{remark}\label{remark_HB_method} Under the hypotheses of Theorem \ref{thm_main}, for $F$ as in (\ref{F_dfn}), $V(F(Y,\X))$ is absolutely irreducible (following similar reasoning to Remark \ref{remark_V_irred}). As a result of this lemma, we can obtain the bound claimed in Theorem \ref{thm_main} as long as all later dependence on $\|F\|$ is at most logarithmic in $\|F\|,$ which we track as the argument proceeds. 

\end{remark}
\begin{proof}
The method of proof follows \cite[Thm. 4]{HB02}, or the recent similar result \cite[Lemma 2.1]{BonBro22}.
Fix $n,D,e \geq  1$.
We start by considering the set of monomials
\[
\mathcal{E}:=\left\{Y^{d_{Y}}X_{1}^{d_{1}}\cdots X_{n}^{d_{n}}:\text{ }d_{Y}e+\sum_{i=1}^{n}d_{i}=De\right\},
\]
in which the degrees $d_Y,d_1,\ldots,d_n$ vary over all non-negative integers satisfying $d_Ye + \sum d_i = De$.
It is easy to see that $|\mathcal{E}|\leq (De)^{n+1}$. 

Let $B\geq 1$ be fixed. Let $\bfv$ denote coordinates $(y,x_{1},\ldots,x_{n})$ and let $\{\bfv_{1},\ldots \bfv_{N}\}$ enumerate the set of points that are solutions to $G(Y,\X)=0$, with each of the last $n$ coordinates of $\bfv_j$ lying in $[-B,B]$. Note that these count each $\X \in [-B,B]^n$ for which $G(Y,\X)$ is solvable at least once, so that $N(G,B) \leq N \leq D \cdot N(G,B).$ (For the upper bound, we recall that the coefficient of $Y^{D}$ in $G(Y,\Xbf)$ is nonzero, so that any given $\X$ can correspond to at most $D$ such $Y$.)
Then, we construct the $N\times |\mathcal{E}|$ matrix
\[
\bfC=(\bfv_{i}^{\bfe})_{\substack{1\leq i\leq N\\ \bfe\in\mathcal{E}}}.
\]
Notice that $\rank \bfC\leq |\mathcal{E}|-1$, since the vector $\bfa\in\mathbb{Z}^{|\mathcal{E}|} \setminus\{0\}$ whose entries correspond to the coefficients of $G(Y,\X)$ is such that $\bfC\bfa=\boldsymbol{0}$.  Moreover, $\bfa$ is primitive since the coefficient associated to $Y^{D}$ is $1$. 
Now the strategy is to find another nonzero vector $\bfb$ in the nullspace of $\mathbf{C}$ and show that if $\bfb$ is in the span of $\abf$ then $\|G\|$ is small, and if $\bfb$ is not in the span of $\abf$ then 
we have an improved count for $N(G,B)$. We may assume henceforward that  $|\mathcal{E}| \leq N,$ since otherwise we already have the upper bound  $N(G,B) \leq N \leq |\mathcal{E}| \leq (De)^{n+1},$ which suffices for the lemma.

 If $\rank \bfC \leq |\mathcal{E}|-2$, then the nullspace has dimension at least 2, and we can take $\bfb \in \Z^{|\mathcal{E}|}$ to be any element in the nullspace that is not in the span of $\bfa$. 
 Let $H(Y,\X)$ be the polynomial defined by the coefficients corresponding to the vector $\bfb$ and consider the polynomial $R(\X)=\text{Res}(G(Y,\X),H(Y,\X))$, which is a polynomial in $\X$ of degree $\ll_{D,e,n} 1$.  (See e.g. \cite[Ch 12]{GKZ08}, which we apply to take the resultant of two polynomials in the variable $Y$, whose coefficients are determined by $\X$.)
We claim that $R(\X)\not\equiv 0$: indeed, if $R(\X)\equiv 0$, then $G$ and $H$ would share an irreducible component.  Since $G(Y,\Xbf)=0$ is irreducible, and $\deg H \leq De  = \deg G$,  it would follow that $G$ is a constant multiple of $H$,  but this is not possible since we are assuming that $\bfa$ and $\bfb$ are not proportional.
Thus $R(\X)\not\equiv 0$. Moreover, observe that for any $\bfx\in\mathbb{Z}^{n}$
\[
R(\bfx)=0\Leftrightarrow G(Y,\bfx)\text{ and }H(Y,\bfx) \text{ have a common root}.
\]
 Note that any $\xbf$ such that $G(y,\xbf)=0$ is solvable contributes  at least one row to the matrix $\bfC$;  each such row also corresponds to a solution to $H(y,\xbf)=0$.
Thus it follows that 
\[
\begin{split}
N(G,B)&=|\{\bfx\in [-B,B]^{n}: \exists y\in\mathbb{Z}\text{ such that } G(y,\bfx)=H(y,\bfx)=0\}|\\&\leq|\{\bfx\in [-B,B]^{n}:\text{ } R(\bfx)=0\}|\\&\ll_{n,D,e,} B^{n-1},
\end{split}
\]
with an implicit constant independent of the coefficients of $R$, via an application of a trivial counting bound for the nonzero polynomial $R$.  (This bound is sometimes called the Schwartz-Zippel bound, and a proof can be found  in \cite[Theorem $1$]{HB02}; we remark that although in that context the polynomial under consideration is absolutely irreducible, the method of proof only requires that it is not identically zero.)

 The remaining case is when $\rank \bfC = |\mathcal{E}|-1$, so that all $|\mathcal{E}| \times |\mathcal{E}|$ minors vanish, but at least one $(|\mathcal{E}|-1) \times (|\mathcal{E}|-1)$ minor does not; we claim there is 
  a nonzero $\bfb \in \Z^{|\mathcal{E}|}$ in the nullspace of $\mathbf{C}$ such that  $|\bfb|=O(B^{De|\mathcal{E}|})=O(B^{(De)^{n+2}})$. 
  If so, then since $\bfa$ is primitive (and $\bfb$ must be proportional to $\bfa$) it follows that $|\bfa|\leq |\bfb|\ll B^{(De)^{n+2}}$. This shows that $\|G\|\ll B^{(De)^{n+2}}$ as claimed.

An appropriate $\bfb$ can be constructed with entries that are $(|\mathcal{E}|-1) \times (|\mathcal{E}|-1)$ minors, so that the size estimate $|\bfb|=O(B^{De|\mathcal{E}|})$ follows from the fact that 
 each entry of $\bfC$ is $O(B^{De}).$ For completeness, we sketch this construction. 
Without loss of generality, we can let $\mathbf{C}'$ denote the top $|\mathcal{E}| \times |\mathcal{E}|$ submatrix in $\mathbf{C}$, and assume that the minor $\mathbf{C}'_{1,1}$ (obtained by omitting the first row and first column of $\mathbf{C}'$) is nonzero. Define a vector $\mathbf{b}$ as follows: for each $1 \leq j \leq |\mathcal{E}|,$ define the entry $b_j$ to be the $(1,j)$-th cofactor of $\mathbf{C'}$; in particular $b_1 \neq 0$ so $\mathbf{b}$ is nonzero, and $|\bfb|=O(B^{De(|\mathcal{E}|-1)})=O(B^{De|\mathcal{E}|})$. We now show that $\bfb$ is in the nullspace of $\mathbf{C}$. Let $\mathbf{r}_i$ denote the $i$-th row of $\Cbf$; then for each $1 \leq i \leq N$, 
\beq\label{rb}\mathbf{r}_i \cdot \bfb= \det \left( \begin{array}{c}
\mathbf{r}_i\\
\mathbf{r}_2 \\
\vdots\\
\mathbf{r}_{|\mathcal{E}|}
\end{array} \right) =0.
\eeq
Indeed, for $i=1$ or $i > |\mathcal{E}|,$ up to sign, $\mathbf{r}_i \cdot \bfb$ is an $|\mathcal{E}| \times |\mathcal{E}|$ minor of $\mathbf{C}$, and all such minors vanish since $\mathrm{rank} \Cbf < |\mathcal{E}|.$ For $2 \leq i \leq |\mathcal{E}|$, the matrix (\ref{rb}) has two identical rows. Thus $\mathbf{C}\bfb=\boldsymbol{0}$.

\end{proof}

\section{Preliminaries on the sieve lemma}
In this section we gather together two preliminary steps: first, we prove the sieve inequality in Lemma \ref{lemma_sieve}; for $m=1$ we provide an alternative proof, conditional on GRH. Second, we formulate an equivalent nonsingularity condition in unweighted projective space.  We also make preliminary remarks on the sieving set.
\subsection{Proof of the polynomial sieve lemma}\label{sec_sieve_proof}

To prove Lemma \ref{lemma_sieve}, observe that
\[
\mathcal{S}(F,B)=\sum_{\bfk:f_{d}(\bfk)=0}W(\bfk)+\sum_{\substack{\bfk\in\mathbb{Z}^{n}:\\f_{d}(\bfk)\neq 0\\ F(y,\bfk)=0 \text{ solvable}}}W(\bfk),
\]
since within the first term, $y=0$ is always a solution to $F(y,\bfk)=0$.
We consider the weighted sum
\beq\label{sieve_lemma_square}
\sum_{\bfk :f_{d}(\bfk)\neq 0}W(\bfk)\left(\sum_{p\in\mathcal{P}}(\nu_{p}(\bfk)-1)\right)^{2}.
\eeq
Fix $\bfk$ such that $f_d(\bfk) \neq 0$ and the polynomial $F(Y,\bfk)$ is solvable over $\mathbb{Z}$, so that  there exists $y_0\in\mathbb{Z}$ such that $F(y_0,\bfk)=0$. For any $p \in \mathcal{P}$ such that $p \ndiv f_d(\bfk)$, then $y_0  \not\equiv 0 \mod p$. Then since $p \equiv 1\mod m$, and due to the structure of $F$ in (\ref{F_dfn}), we have that $\{y_0,\gamma_{p}y_0,\ldots,\gamma_{p}^{m-1}y_0\}$ are distinct solutions of $F(Y,\bfk)\con 0 \modd{p}$, where $\gamma_{p}^{m}\equiv 1\mod p$ and $\ga_p$ is a primitive $m$-th root of unity in $\F_p$. In particular, for such $p$, $\nu_p(\bfk) \geq m$. 
Consequently, for each $\bfk $ such that $f_{d}(\bfk)\neq 0$ and $F(Y,\bfk)$ is solvable, we have that
\beq\label{Pm}
\sum_{p\in\mathcal{P}}(\nu_{p}(\bfk)-1)\geq (m-1)\sum_{p\in\mathcal P, p\nmid f_{d}(\bfk)}1\gg_m P- \sum_{p\in\mathcal P, p \mid f_{d}(\bfk)} 1 \geq (1/2) P,
\eeq
as long as $P \gg_{m,e,d} \max\{ \log \|f_d\|, \log B\}$.
 The last step follows since the number $\omega(f_d(\bfk))$ of distinct prime divisors of $f_d(\bfk) \neq 0$ is at most 
\begin{align*}
\omega(f_d(\bfk))&\ll \log (f_d(\bfk))/ \log \log (f_d(\bfk))\\
&\ll \log (\|f_d\| B^{dem})   \\ &\ll_{m,e,d} \log \|f_d\|  + \log B.
\end{align*}
 Thus the last inequality in (\ref{Pm}) holds as long as \beq\label{P_big}
 P \gg_{m,e,d} \max\{ \log \|f_d\|, \log B\},
 \eeq
 leading to the corresponding hypothesis in the lemma.

From  (\ref{Pm}) and the non-negativity of the weight $W$, we see that
\[
P^{2}\sum_{\substack{\bfk\in\mathbb{Z}^{n}:\\f_{d}(\bfk)\neq 0\\ F(y,\bfk)=0 \text{ solvable}}}W(\bfk)
	\ll\sum_{\bfk :f_{d}(\bfk)\neq 0}W(\bfk)\left(\sum_{p\in\mathcal{P}}(\nu_{p}(\bfk)-1)\right)^{2}
	\leq \sum_{\bfk }W(\bfk)\left(\sum_{p\in\mathcal{P}}(\nu_{p}(\bfk)-1)\right)^{2}.
\]
Opening the square on the right-hand side, the contribution from $p=q \in \mathcal{P}$ is 
\[ \sum_{p \in \Pcal} \sum_{\bfk  }W(\bfk) (\nu_{p}(\bfk)-1)^{2} \ll_{m,d} P \sum_{\bfk  }W(\bfk)  ,
\]
since $\nu_p(\bfk) \leq md$ for all $\bfk$, as previously mentioned.
The contribution from $p \neq q \in \Pcal$ is bounded in absolute value by
\[ \sum_{p \neq q \in \Pcal} |\sum_{\bfk }W(\bfk) (\nu_{p}(\bfk)-1) (\nu_{q}(\bfk)-1) |.
\]
Assembling all these terms, we see that Lemma \ref{lemma_sieve} is proved.

 \begin{remark}
When we apply Lemma \ref{lemma_sieve} to prove Theorem \ref{thm_main}, we can assume that  $\|f_d \| \leq \|F\|\ll  B^{(mde)^{n+2}}$, by Lemma \ref{lemma_F_small}. This will allow us to verify that (\ref{P_big}) holds for our choice of sieving set, as we will verify in \S \ref{sec_concluding} when we choose $Q$ in (\ref{Q_choice}). 
\end{remark}

\subsection{Alternative proof when $m=1$, conditional on GRH}\label{sec_GRH}

Recall from \S \ref{sec_intro_poly_sep} the general problem of counting $\xbf \in [-B,B]^n$ such that $G(y,\xbf)=0$ is solvable in $\Z$, with $G(Y,\Xbf)$ of degree $D$ as in (\ref{G_dfn_intro}). In our main work in this paper, we assume that   $D=md$ with $m \geq 2$, and $G$ is a polynomial in $Y^m$.   This  additional structure allowed us to choose a sieving set $\Pcal \subset [Q,2Q]$ of primes $p \con 1 \modd{m}$, so that all the $m$-th roots of unity are present in $\F_p$, for each $p \in \Pcal$. With this property, we could define sieve weights that exhibit an appropriate lower bound in the form (\ref{Pm}) for most $\kbf$ in the support of $W(\kbf)$ and a positive proportion of primes.  

Nevertheless, we can proceed by a different argument to  develop a sieve lemma to bound the number of $\x \in [-B,B]^n$ such that  $G(y,\xbf)=0$ is solvable over $\Z$, with no condition on the degree $D$; that is, to prove a  version of Lemma \ref{lemma_sieve} in the case $m=1$. As a first step, we naturally try to introduce a system of weights, according to a fixed set of primes. Let us take $\mathcal{P}= \{Q\leq p\leq 2Q: p \text{ prime}\}$ for some parameter $Q$ to be chosen optimally with respect to $B$. In particular, by the prime number theorem, $|\Pcal| \gg Q (\log Q)^{-1}$ for all $Q \gg 1$. Fix $\bfk \in \Z^n$. 
For each prime $p \in \mathcal{P}$, set 
\[\nu_{p}(\bfk)=|\{y\in\mathbb{F}_{p}: G(y,\bfk)=0 \modd{p}\}|.\] 
Since $G(y,\bfk)$ contains the term $y^D,$ it is not the zero polynomial in $y$, and $\nu_p(\kbf) \leq D$.
Consider, as in the proof of Lemma \ref{lemma_sieve} above, the weighted sum
\begin{equation}\label{sieve_weight}
\sum_{\bfk: f_D(\kbf) \neq 0}W(\bfk)\left(\sum_{p\in\mathcal{P}}(\nu_{p}(\bfk)-1)\right)^{2}.
\end{equation}
In order to deduce a sieve lemma, we need a lower bound for the arithmetic weight (the squared term), for those $\bfk$ for which $f_D(\bfk) \neq 0$ and $G(Y,\bfk)=0$ is solvable over $\Z$.

Here is one approach. Let $\kbf$ be fixed, with $f_D(\bfk) \neq 0$ and $G(Y,\bfk)=0$  solvable over $\mathbb{Z}$, and $\kbf$ in the support of $W$. Then $G(Y,\bfk)=(Y-y_0)\tilde{g}_{\bfk}(Y)$ for some $y_0\in\mathbb{Z} \setminus \{0\}$ and some (monic) $\tilde{g}_{\bfk}(Y) \in \mathbb{Z}[Y]$ of degree $D-1$. 
For each such $\kbf$, we can obtain a suitable lower bound for the arithmetic weight in (\ref{sieve_weight}) as long as for a positive proportion of $p \in \Pcal$, $\tilde{g}_\kbf$ has a root over $\F_p$. Let $g_{\bfk}$ be an irreducible factor of $\tilde{g}_\kbf$. 
Let $F_{\bfk}$ denote the splitting field of $g_{\bfk}$ over $\Q$, say $F_{\bfk} = \Q(\al_{\bfk}).$ Since $g_{\bfk}$ is irreducible, then it is the minimal polynomial of $\al_{\bfk}$ in $\Z[Y]$, and it is separable (since we are working over characteristic zero), and the splitting field is Galois over $\Q$. 
By Dedekind's theorem, for all $p \ndiv [\mathcal{O}_{F_{\bfk}} : \Z[\al_\kbf]]$, $g_\kbf$ splits completely over $\F_p$ precisely when $(p) = p\mathcal{O}_{F_{\bfk}}$ splits completely in $F_{\bfk}$; see e.g. \cite[Thm. 27 p. 79]{Mar77}.  
Then 
\[\sum_{p\in\mathcal{P}}(\nu_{p}(\bfk)-1)
    =\sum_{p\in\mathcal{P}}|\{y\in\mathbb{F}_{p}:\tilde{g}_{\bfk}(y)=0\}|
 \geq \sum_{p\in\mathcal{P}}|\{y\in\mathbb{F}_{p}:g_{\bfk}(y)=0\}|.
 \]
If $g_{\bfk}$ is linear in $\Z[Y]$, this sum is of size $|\mathcal{P}|$, which suffices. If $\deg g_{\bfk} \geq 2$, we continue to argue that
\begin{align}\sum_{p\in\mathcal{P}}(\nu_{p}(\bfk)-1) 
 & \geq \deg(g_{\kbf})|\{p\in\mathcal{P}: g_{\bfk}(Y) \text{ completely split over $\mathbb{F}_{p}$}\}|\nonumber \\
  & \geq |\{ p \in \mathcal{P}: \text{$ p\mathcal{O}_{F_{\bfk}}$ splits completely in $F_{\bfk}$}\}| - |\{ p \in \mathcal{P}: p| [\mathcal{O}_{F_{\bfk}} : \Z[\al_\kbf]]\}|. \label{arith_lower}
\end{align}
Let 
\[\pi_{\kbf}(Q)=|\{ p \leq Q: \text{$ p\mathcal{O}_{F_{\bfk}}$ splits completely in $F_{\bfk}$}\}|\]
and $N(\kbf) = |\{  p| [\mathcal{O}_{F_{\bfk}} : \Z[\al_\kbf]]\}|$.
The Chebotarev density theorem, in the unconditional form of \cite[Thm. 1.3]{LagOdl77}, shows that 
\beq\label{Chebotarev} \left|\pi_{\kbf}(Q)     - \frac{1}{|G_{\bfk}|}\frac{Q}{\log Q} \right| =  \frac{1}{|G_{\bfk}|}\frac{Q^{\be_0}}{\log Q^{\be_0}}  + O_{D,A}(Q (\log Q)^{-A}) 
\eeq
     for every $A \geq 2$, as long as $Q \geq \exp (10 \deg F_{\bfk} (\log |D(F_{\bfk})|)^2).$ Here
   $G_{\bfk}$ is the Galois group  $\mathrm{Gal}(F_{\bfk}/\Q)$, $D(F_{\bfk})$ is the discriminant of the splitting field $F_{\bfk}/\Q,$ and  $\deg F_{\kbf}= \deg|F_\kbf/\Q|$ is the degree of the extension. The implicit constant in the error term depends only on $A$ and $\deg F_{\kbf} =|G_{\bfk}|\leq (D-1)!$. The real number $1/2<\beta_0<1$, if it exists, is the (real, simple) exceptional zero of the associated Dedekind zeta function $\zeta_{F_{\bfk}};$ if no exceptional zero exists, that term does not appear in the result. 
   
   In particular, under the assumption of GRH for $\zeta_{F_{\bfk}},$ Lagarias and Odlyzko's Theorem 1.1 in \cite{LagOdl77} (in the refined form of Serre \cite[Thm. 4]{Ser81}) shows that for any $Q>2$, the entire right-hand side of (\ref{Chebotarev}) may be replaced by 
   \[O(|G_{\bfk}|^{-1} Q^{1/2} \log (|D(F_{\bfk})| Q^{\deg F_{\bfk}})) = O_D(Q^{1/2} \log  Q) + O_D(Q^{1/2} \log |D(F_{\bfk})|),\]
   in which the implied constant is absolute and effectively computable. There exists a constant $Q_0(D)$ depending only on $D$ such that the first term is $\leq \frac{1}{4}\frac{1}{(D-1)!} Q(\log Q)^{-1}$ for all $Q \geq Q_0(D).$ The second term is also $\leq \frac{1}{4}\frac{1}{(D-1)!} Q(\log Q)^{-1}$ if for example $Q \geq Q_1(D) (\log D(F_\kbf))^{\al_0}$ for a constant $Q_1(D)$ and some fixed $\al_0>2.$
   This shows that under GRH, for all $Q \gg_D (\log D(F_\kbf))^{\al_0}$ some fixed $\al_0>2$, 
   \beq\label{pi_lower}
   \pi_\kbf(Q) - \pi_\kbf(Q/2) \gg_D Q/\log Q \gg_D|\Pcal|.\eeq

 Two tasks remain in order to complete a lower bound for (\ref{arith_lower}): (i) to bound $D(F_\kbf)$ from above, so that the lower bound $Q \gg_D (\log D(F_\kbf))^{\al_0}$  can be made uniform over $\kbf$, and (ii)  to count 
\[N(\kbf) = |\{  p| [\mathcal{O}_{F_{\bfk}} : \Z[\al_\kbf]]\}| = \om([\mathcal{O}_{F_{\bfk}} : \Z[\al_\kbf]]) \ll \log [\mathcal{O}_{F_{\bfk}}  : \Z[\al_\kbf]]/ \log \log [\mathcal{O}_{F_{\bfk}} : \Z[\al_\kbf]].\]
We note the relation    \beq\label{D_relation} D(F_\kbf) [\mathcal{O}_{F_{\bfk}} : \Z[\al_\kbf]]^2 =\mathrm{Disc}(g_\kbf),
 \eeq
which holds by \cite[Remark 2.25 and Eqn. (8) on p. 38]{Mil20x}. 
(Since $g_\kbf$ was assumed to be irreducible and we are in characteristic zero, then $g_\kbf$ is separable and $\mathrm{Disc}(g_\kbf) \neq 0.$)
 Thus for both remaining tasks, it suffices to bound $\mathrm{Disc}(g_\kbf)$ from above, since by (\ref{D_relation}) both  
 \[ N(\kbf) \ll \log \Disc(g_\kbf), \qquad \log D(F_\kbf) \leq \log \Disc(g_\kbf).\]

Now $\Disc(g_\kbf)$ (the resultant of $g_\kbf(Y)$ and $g_\kbf'(Y)$, as defined in \cite[Prop. 1.1, Ch. 13]{GKZ08}) is a polynomial in the coefficients of $g_{\kbf}$  with degree bounded in terms of $D$. The coefficients of $g_\kbf$ are polynomials in   $\kbf$ and the coefficients of $G(Y,\Xbf)$ with degree at most $D$. Since we only consider $\kbf$ in the support of $W$, $|\kbf| \ll B$, and the coefficients of $g_{\kbf}$ are $\ll \|G\|B^{D}.$
Thus 
\[ \log \Disc (g_\kbf) \ll_D \log \|G\|  + \log B.\]
  
In combination with (\ref{pi_lower}), we can conclude in (\ref{arith_lower}) that for some constant $C_D$,
 \[ \sum_{p \in \mathcal{P}}(\nu_p(\kbf)-1) \gg_D Q/\log Q - C_D (\log \|G\| + \log B),\]
for all $Q \geq C'_D \max\{(\log \|G\|)^{\al_0}, (\log B)^{\al_0}\}$ for some $\al_0>2$. By  taking $C'_D$ sufficiently large, we achieve
$ \sum_{p \in \mathcal{P}}(\nu_p(\kbf)-1) \gg |\mathcal{P}|=P.$
This shows that conditional on GRH, 
\[
P^{2}\sum_{\substack{\bfk\in\mathbb{Z}^{n}:\\f_{D}(\bfk)\neq 0\\ G(y,\bfk) =0\text{ solvable}}}W(\bfk)
	\ll\sum_{\bfk :f_{D}(\bfk)\neq 0}W(\bfk)\left(\sum_{p\in\mathcal{P}}(\nu_{p}(\bfk)-1)\right)^{2}
	\leq \sum_{\bfk }W(\bfk)\left(\sum_{p\in\mathcal{P}}(\nu_{p}(\bfk)-1)\right)^{2}.
\]
 From here, the remainder of the proof used above for Lemma \ref{lemma_sieve} can be repeated, and this completes the proof of the claim in Remark \ref{remark_GRH}.

\subsection{Associated variety in unweighted projective space}\label{sec_smooth}

It is a hypothesis of Theorem \ref{thm_main} that the weighted hypersurface $V(F(Y,\X)) \subset \mathbb{P}(e,1,\ldots,1)$, defined by $F(Y,\Xbf)=0$, is nonsingular over $\C$. It is convenient to relate $V(F(Y,\X))$ to a variety in unweighted projective space. 
We claim that for 
\[
F(Y,\X)=Y^{dm}+Y^{(d-1)m}f_{1}(\X)+\ldots + f_{d}(\X),
\]
then $V(F(Y,\X))\subset\mathbb{P}(e,1,\ldots,1)$ is nonsingular  if and only if $V(F(Z^{e},\X))\subset\mathbb{P}^{n}$ is nonsingular. Here, we again apply the assumption $m \geq 2$. Indeed the weighted projective variety is nonsingular if and only if the only solution of
\begin{equation}
\begin{cases}
F(Y,\X)=0\\
\frac{\partial F}{\partial Y}(Y,\X)= \sum_{i =0}^{d-1}f_{i}(\X)\cdot m(d-i)Y^{m(d-i)-1}=0\\
\frac{\partial F}{\partial X_{1}}(Y,\X)=0\\
\vdots\\
\frac{\partial F}{\partial X_{n}}(Y,\X)=0
\end{cases}
\label{eq : jacY}
\end{equation}
on $\mathbb{A}^{n+1}$ is the point $P=\boldsymbol{0}$.
(By convention we set $f_0(\Xbf)=1.$) Similarly, the projective variety $V(F(Z^{e},\X))$ is nonsingular if and only if
the only solution of
\begin{equation} 
\begin{cases}
F(Z^{e},\X)=0\\
\frac{\partial F}{\partial Z}(Z^e,\X)= \sum_{i =0}^{d-1}f_{i}(\X)\cdot me(d-i)Z^{em(d-i)-1}=0\\
\frac{\partial F}{\partial X_{1}}(Z^{e},\X)=0\\
\vdots\\
\frac{\partial F}{\partial X_{n}}(Z^{e},\X)=0
\end{cases}
\label{eq : jacZ}
\end{equation}
on $\mathbb{A}^{n+1}$ is the point $P=\boldsymbol{0}$. Moreover, note that
\begin{align}
\frac{\partial F}{\partial Y}(Y,\X)&=mY^{m-1}\sum_{i=0}^{d-1}f_{i}(\X)(d-i)Y^{m(d-i-1)}\label{Y_factor}\\
\frac{\partial F}{\partial Z}(Z^{e},\X)&=emZ^{em-1}\sum_{i=0}^{d-1}f_{i}(\X)(d-i)Z^{em(d-i-1)}.\nonumber
\end{align}
We will momentarily use this to confirm that if $m\geq 2$,  a nonzero solution (say $P=(y,\bfx) \in \mathbb{A}^{n+1}$) to  $(\ref{eq : jacY})$ exists if and only if   a solution (namely $Q=(y^{1/e},\x) \in \mathbb{A}^{n+1}$) to $(\ref{eq : jacZ})$ exists.

To clarify the role of the assumption $m \geq 2$, let us briefly make a general observation.
In general, let a polynomial  $G(Y,\X)$ be given as in (\ref{G_dfn_intro}) and assume $V(G(Y,\X))\subset \mathbb{P}(e,1,\ldots,1)$ is nonsingular; we may assume $e \geq 2$ (since otherwise the variety is already unweighted).
Then we claim $V(G(Z^e,\X))$ is nonsingular (as a projective variety) if and only if $V(G(Y,\X))\cap V(Y)$ is nonsingular (as a weighted projective variety). 
By the chain rule, 
\[\frac{\partial G}{\partial Z}(Z^{e},\X)=eZ^{e-1}(\frac{\partial G}{\partial Y})(Z^{e},\X).\]
Observe that
\begin{align} \mathrm{Sing}(V(G(Z^e,\X)))
&= \{(z,\xbf)\in \mathbb{P}^n : \nabla_{Z,\Xbf}G(z^e,\xbf)=\zerobf\} \nonumber \\
&= 
\{ (0,\xbf)\in \mathbb{P}^n : \nabla_\Xbf G(0,\xbf)=\zerobf\} \cup \{(z,\xbf)\in \mathbb{P}^n : \nabla_{Y,\Xbf}G(z^e,\xbf)=\zerobf\} \label{two_sets}\\
&= 
\{ (0,\xbf)\in \mathbb{P}^n : \nabla_\Xbf G(0,\xbf)=\zerobf\} \cup \emptyset \nonumber
\end{align}
under the assumption that $V(G(Y,\X))$ is nonsingular.
On the other hand, by the Jacobian criterion,
\[ \mathrm{Sing}(V(G(Y,\X))\cap V(Y)) = \{ (0,\xbf) \in \mathbb{P}^n :  \nabla_{\Xbf}G(0,\xbf)=\zerobf\}.\]
(Here we have used that $G(0,\Xbf)$ is itself homogeneous in $\Xbf$, so that $\nabla_XG(0,\Xbf)=0$ implies $G(0,\Xbf)=0$ by Euler's identity.)
Since the singular sets are identical, this proves the claim. 

Let us apply this in our case with $G$ taken to be the polynomial $F(Y,\Xbf)$, with $V(F(Y,\Xbf))$ assumed to be nonsingular. We consider whether there are any $(0,\xbf) \in \mathbb{P}^n$ such that $\nabla_\Xbf F(0,\xbf)=0.$ Supposing such $(0,\xbf)$ exists, it must be the case that $(\frac{\partial F}{\partial Y})(0,\xbf) \neq 0,$ since otherwise $(0,\xbf)$ would be a singular point  on $V(F(Y,\Xbf)).$
If $m \geq 2$, then due to the leading factor $Y^{m-1}$ in (\ref{Y_factor}), any point $(0,\xbf) \in \mathbb{P}^n$ must lead to
$(\frac{\partial F}{\partial Y})(0,\xbf) =0$. Consequently there can be no such $(0,\xbf$), and $\mathrm{Sing}(V(F(Y,\X))\cap V(Y))$ must be empty. Hence by the general argument above, so is $\mathrm{Sing}(V(F(Z^e,\Xbf))$. 
In conclusion, if $m \geq 2,$ $V(F(Y,\Xbf))$ being nonsingular implies $V(F(Z^e,\Xbf))$ is nonsingular.

However if $m=1$, there is no leading factor of $Y$ in (\ref{Y_factor}), and indeed at $(0,\xbf)$, (\ref{Y_factor}) evaluates to $f_{d-1}(\xbf)$.
Thus points $(0,\xbf)$ for which $f_{d-1}(\xbf)\neq 0$ and $\nabla_{\Xbf} F(0,\xbf)=0$ can lead to singular points on $V(F(Y,\X))\cap V(Y)$ and hence to singular points on $F(F(Z^e,\Xbf))$. (Nevertheless, there cannot be too many singular points, as we will observe in (\ref{sing_dim}) below that the singular locus has at most dimension 0.)

In the other direction, suppose that $V(F(Z^e,\Xbf))$ is nonsingular, so that as computed in (\ref{two_sets}), 
\[\mathrm{Sing}(V(F(Z^e,\X)))
= 
\{ (0,\xbf)\in \mathbb{P}^n : \nabla_\Xbf F(0,\xbf)=\zerobf\} \cup \{(z,\xbf)\in \mathbb{P}^n : \nabla_{Y,\Xbf}F(z^e,\xbf)=\zerobf\}\]
is empty. If there were a point $(y,\xbf)$ in $\mathrm{Sing}(V(Y,\Xbf))$ then if $y =0$ this would produce an element in the first set on the right-hand side, while if $y \neq 0$ then taking $z=y^{1/e}$ (working over $\C$) would produce a point in the second set on the right-hand side. Thus $V(F(Y,\Xbf))$ must be nonsingular (and here we did not need to apply $m \geq 2$).

\begin{remark}\label{remark_d1}
In the special case that $d=1$, then   $F(Y,\Xbf)=Y^{m} + f_1(\Xbf).$ Thus $V(F(Y,\Xbf)) \subset \mathbb{P}(e,1,\ldots,1)$ is nonsingular if and only if $V(Z^{em} +f_1(\Xbf)) \subset \mathbb{P}^n$ is nonsingular, with $f_1 \not\con 0$ homogeneous of degree $em.$ This occurs if and only if $V(f_1(\Xbf))\subset \mathbb{P}^{n-1}$ is nonsingular; in this special case, the problem we consider falls in the scope of the work in \cite[Theorem 1.1]{Bon21}, which proves this  case of Theorem \ref{thm_main}. Our method of proof works regardless, so we allow $d=1$ as we continue. 
\end{remark}

\begin{remark}\label{remark_V_irred}

Recall the affine hypersurface $\mathcal{V} \subset \mathbb{A}_\C^{n+1}$ defined in (\ref{Vcal}) according to the polynomial $F(Y,\X)$. We note that $\mathcal{V}$ is irreducible under the conditions of Theorem \ref{thm_main}. Suppose it is reducible, so that $F(Y,\X)=G(Y,\X)H(Y,\X)$ for some nonconstant polynomials. Then $F(Z^e,\X)=G(Z^e,\X)H(Z^e,\X)$ so that the projective variety $V(F(Z^e,\X))$ is reducible. Consequently, by \cite[Lemma 11.1]{BCLP23}, $V(F(Z^e,\X))$ is singular, which is a contradiction because by the discussion above, $V(F(Y,\X))$ is nonsingular if and only if $V(F(Z^e,\X))$ is nonsingular. 
\end{remark}

\subsection{Initial considerations of the sieving set}\label{sec_sieve_first}
We suppose that $Q=B^\kappa$ for some $0<\kappa \leq 1$ to be chosen later (see (\ref{Q_choice})). 
We will choose a sieving set 
\[ \Pcal \subset  [Q,2Q] \]
 comprised of primes with certain properties. 
 In the special case that $(e,m)=1$, it is sensible to restrict our attention to a set $\Pcal_0$ of primes in $[Q,2Q]$ such that:\\
(i) $p \con 1 \modd{m}$ (recalling $m \geq 2$) and \\
(ii)  $p\equiv 2\mod e$, and 
\\
(iii) the reduction of $V(F(Y,\X))$ as a weighted variety over $\overline{\F}_p$ is nonsingular.  

 The first criterion (i) we have used in the proof of the sieve lemma (Lemma \ref{lemma_sieve}). The second criterion (ii) ensures that $(e,p-1)=1$ so that every $y \in \F_p$ satisfies $y=z^e$ for some $z \in \F_p$.  Then for each $p \in \Pcal$, we can simply consider the reduction $V(F(Z^e,\Xbf)) \subset \mathbb{P}_{\F_p}^n$ in place of the weighted variety, so that (iii) is equivalent to:\\
 (iii') the reduction of $V(F(Z^e,\Xbf)) \subset \mathbb{P}_{\overline{\F}_p}^n$ is nonsingular.
 
By the Chinese remainder theorem and the Siegel--Walfisz theorem on primes in arithmetic progressions, under the assumption that $(e,m)=1$, there are $\gg_{m,e} Q/\log Q$ primes that satisfy (i) and (ii) in any dyadic region $[Q,2Q],$ for all $Q$ sufficiently large. We could then choose the sieving set $\Pcal_0$ to be the subset of such primes for which (iii') holds; the remaining task is to show there are sufficiently few primes that violate (iii'). 

 Recall from \S \ref{sec_smooth} that $V(F(Y,\Xbf))$ is nonsingular over $\C$ (as a weighted projective variety) if and only if $V(F(Z^e,\Xbf)) \subset \mathbb{P}^n$ is nonsingular over $\C$. Thus under the hypothesis of Theorem \ref{thm_main}, the latter is nonsingular, and consequently 
there are no nontrivial simultaneous solutions of the system (\ref{eq : jacZ}), and thus the resultant \[
r:=\mathrm{Res}(F, \frac{\partial F}{\partial Z}, \frac{\partial F}{\partial X_{1}}, \ldots, \frac{\partial F}{\partial X_{n}})
\]
of those $n+2$ polynomials in $n+1$ variables is a nonzero integer.
 Moreover, by \cite[Prop. 1.1, Ch. 13]{GKZ08}, $r$ is a polynomial in the coefficients of $F$ with degree bounded in terms of $m, e, d$. By \cite[Section IV]{Cha93}, the reduction $V_p(F(Z^e,\X))$ of $V(F(Z^e,\X))$ modulo $p$ is singular precisely when $p|r$, which can only occur for at most $\omega (r)$ primes,
where
\beq\label{few_primes}
\om(r) \ll \log r/\log \log r \ll_{m,e,d}\log \|F\|.
\eeq
(Notice that the argument in this paragraph made no assumption on the relative primality of $e$ and $m$.)

In particular, if $(e,m)=1$, then as long as  $Q$ is sufficiently large, say $Q \gg_{m,e,d} (\log \|F\|)^{1+\del_0}$ for any fixed $\del_0>0$ or even $Q \gg_{m,e,d} (\log \|F\|)(\log \log \|F\|)$, we can conclude that $|\Pcal_0| \gg_{m,e,d} Q/\log Q.$ After we choose $Q$ to be a certain power of $B$ (see (\ref{Q_choice})), this will only require a lower bound on $B$ that is on the order of a power of $\log \|F\|$, which we will see can be accommodated by the bound on the right-hand side of our claim in Theorem \ref{thm_main}.

These remarks all apply in the case that $(e,m)=1$. However, we can also argue more generally without this assumption, as we demonstrate in the next section, by working not with $V(F(Z^e,\X))$ as above, but with 
a finite collection of varieties 
$W_i$, defined according to $F(\ga^i z^e,\Xbf)=0$ in $\F_p$, for a certain primitive root $\ga \in \F_p^\times$ (see Lemma \ref{lemma_split_by_root}).
Thus we postpone our definition of the sieving set, in general, until the end of the next section.

\section{Estimates for exponential sums}
In this section we apply the Weil bound to prove an upper bound for the exponential sum $g(\ubf,p)$ (see $(\ref{eq : expsumg})$) in the case that $\ubf$ is each of three types: type zero, good, or bad modulo $p$ (Definition \ref{dfn_types}). 
At the end, in \S \ref{sec_sieve_set} we then define the sieving set $\Pcal$. 

We note the multiplicativity condition
\[
g(\bfu,pq):=\sum_{\bfa \mod pq}(\nu_{p}(\bfa)-1)(\nu_{q}(\bfa)-1)e_{pq}(\langle\bfa,\bfu\rangle)=g(\overline{q}\bfu,p)g(\overline{p}\bfu,q),
\]
where $q\overline{q}\equiv 1\mod p$, and $p\overline{p}\equiv 1\mod q$.
This leads us to study the key exponential sums with prime modulus:
\[
g(\bfu, p):=\sum_{\bfa \in\mathbb{F}_{p}^{n}}(\nu_{p}(\bfa)-1)e_{p}(\langle\bfa,\bfu\rangle).
\]

Let $p$ be a fixed prime of good reduction for $F(Z^e,\Xbf)$, so that  $V(F(Z^{e},\X))\subset\mathbb{P}_{\overline{\mathbb{F}}_{p}}^{n}$ is a nonsingular projective hypersurface. For any point $P \in V(F(Z^{e},\X))$, let $T_P \subseteq \mathbb{P}_{\overline{\mathbb{F}}_{p}}^{n}$ denote the projective tangent space to  $V(F(Z^{e},\X))$ at $P$. A linear space $L$ is tangent to $V(F(Z^{e},\X))$ at $P$ if $T_P \subseteq L$; if $L$ is a hyperplane, this is equivalent to $P$ being a singular point of   $V(F(Z^{e},\X)) \intersect L$ (see \cite[p. 57]{FulLaz81}).

Given $\ubf \in \Z^n$ with $\ubf \not\con \zerobf \modd{p}$, if $V(\langle \X,\bfu\rangle)\subset\mathbb{P}_{\overline{\mathbb{F}}_{p}}^{n}$ is   not tangent to $V(F(Z^{e},\X))$ at  any point (i.e. they intersect  transversely), we simply say $V(\langle \X,\bfu\rangle)$ is  not tangent to $V(F(Z^{e},\X))$; otherwise, we will say they are tangent (and as we will discuss below in (\ref{sing_dim}), there are at most finitely many points at which they are tangent).

Using this terminology, we will  classify $\ubf \in \Z^n$ in terms of three cases: 
\begin{dfn}\label{dfn_types}
For $\bfu\in\mathbb{Z}^{n}$ and $p\in\mathcal{P}$ we say that:
\begin{enumerate}[(i)]
\item $\bfu$ is of type zero mod $p$ if $\bfu \con \boldsymbol{0} \modd{p}$,
\item $\bfu$ is good mod $p$ if $\bfu \not\con \boldsymbol{0} \modd{p}$ and $V(\langle \X,\bfu\rangle)\subset\mathbb{P}_{\overline{\mathbb{F}}_{p}}^{n}$ is not tangent to $V(F(Z^{e},\X))\subset\mathbb{P}_{\overline{\mathbb{F}}_{p}}^{n}$,
\item $\bfu$ is bad mod $p$ if $\bfu \not\con \boldsymbol{0} \modd{p}$, and $V(\langle \X,\bfu\rangle)\subset\mathbb{P}_{\overline{\mathbb{F}}_{p}}^{n}$ is tangent to $V(F(Z^{e},\X))\subset\mathbb{P}_{\overline{\mathbb{F}}_{p}}^{n}$.
\end{enumerate}
\end{dfn}
(The fact that we define these types in relation to $V(F(Z^e,\Xbf))$,   is justified by Lemma \ref{lemma_isom}, below.)
The main result of this section is the following:
\begin{prop}\label{prop_g_sum}
Assume that $p>2$ is a prime of good reduction for $F(Z^e,\Xbf)$, that is $V(F(Z^e,\X)) \subset \mathbb{P}^n_{\overline{\F}_p}$ is nonsingular.
\begin{enumerate}[(i)]
\item If $\ubf$ is type zero modulo $p$ then $g(\ubf,p)\ll p^{n-1/2}$;
\item If $\ubf$ is good modulo $p$ then $g(\ubf,p)\ll p^{n/2}$;
\item If $\ubf$ is bad modulo $p$ then $g(\ubf,p)\ll p^{(n+1)/2}$.
\end{enumerate}
 The implied constants can depend on $n,m,e,d,$ but are independent of $\|F\|,\ubf,p$.
\end{prop}

In a final step  of the proof, we will apply the property that if $V(F(Z^e,\X)) \subset \mathbb{P}^n$ is nonsingular,  any hyperplane $L$ has 
\beq\label{sing_dim} \dim \{ P \in V(F(Z^{e},\X)): T_P \subseteq L\} = \dim(\sing (V(F(Z^{e},\X))\cap L))\leq  0. \eeq 
  Here, by $\dim (\sing (V))$ we mean the dimension of the singular locus of a variety $V \subset \mathbb{P}^n.$
We will apply this in (\ref{Katz_app}) over $\overline{\F}_p$ for $p$ a prime of good reduction for $F(Z^e,\X).$
The result (\ref{sing_dim}) is a special case of Zak's theorem on tangencies as in \cite[Thm. 7.1, Rem. 7.5]{FulLaz81}, valid over any algebraically closed field, or \cite[Lemma 3]{Kat99}, valid over any perfect field. More simply, in our setting (\ref{sing_dim}) can be shown directly, and we do so in Remark \ref{remark_sing_dim}.

As preparation for proving Proposition \ref{prop_g_sum}, we transform $g(\ubf,p)$ into an exponential sum over solutions to $F(y,\abf)=0$ by writing
\[
\begin{split}
g(\bfu,p)&=\sum_{\bfa \in\mathbb{F}_{p}^{n}}\nu_{p}(\bfa)e_{p}(\langle\bfa,\bfu\rangle)-\sum_{\bfa \in\mathbb{F}_{p}^{n}}e_{p}(\langle\bfa,\bfu\rangle)\\&=-\delta_{\bfu =\boldsymbol{0}}\cdot p^{n}+\sum_{\bfa \in\mathbb{F}_{p}^{n}}e_{p}(\langle\bfa,\bfu\rangle)\sum_{\substack{y\in\mathbb{F}_{p}\\F(y,\bfa)=0}}1\\&=-\delta_{\bfu =\boldsymbol{0}}\cdot p^{n}+\sum_{\substack{(y,\bfa ) \in\mathbb{F}_{p}^{n+1}\\ F(y,\bfa)=0}}e_{p}(\langle\bfa,\bfu\rangle),
\end{split}
\]
where $\delta_{\bfu =\boldsymbol{0}}=1$ if  $\bfu\con \boldsymbol{0}\modd{p}$ and   is $0$ otherwise. The  task now is to estimate the sum
\[
g(\ubf,p)+\delta_{\bfu =\boldsymbol{0}}\cdot p^{n}=\sum_{\substack{(y,\bfa ) \in\mathbb{F}_{p}^{n+1}\\ F(y,\bfa)=0}}e_{p}(\langle\bfa,\bfu\rangle).
\]
A barrier to doing this efficiently is that the polynomial $F(Y,\X)$ is not homogeneous (see Remark \ref{remark_F_hom}).  
Recall the definition of $F(Y,\X)$ in (\ref{F_dfn}), and recall the integer $e \geq 1$ fixed in that definition. 
 As a first step, we prove:
\begin{lemma}\label{lemma_split_by_root}
Fix a prime $p>2$. Let $f=(e,p-1)$, and let $\gamma\in\mathbb{F}_{p}^{\times}$ be a primitive $f$-th root of unity. Then
\[
\sum_{\substack{(y,\bfa )\in W}}e_{p}(\langle\bfa,\bfu\rangle)=\frac{1}{f}\sum_{i=0}^{f-1}\sum_{\substack{(z,\bfa ) \in W_{i}}}e_{p}(\langle\bfa,\bfu\rangle),
\]
where
\[
\begin{split}
&W=\{(y,\bfa)\in\mathbb{F}_{p}^{n+1}: F(y,\bfa)=0\}\\
& W_{i}=\{(z,\bfa)\in\mathbb{F}_{p}^{n+1}: F(\gamma^{i}z^{e},\bfa)=0\}, \qquad \text{for $i=0,\ldots, f-1.$}
\end{split}
\]
\label{lem : homogen}
\end{lemma}
(This lemma replaces the  remarks in \S \ref{sec_sieve_first}  that applied in the special case $(e,p-1)=1$.)
\begin{proof}
We start by claiming that for any $y\in\mathbb{F}_{p}^{\times}$ there exists an unique $i\in\{0,\ldots,f-1\}$ and some $z\in\mathbb{F}_{p}^{\times}$ such that $y=\gamma^{i}z^{e}$: we write $e=\ell k$ where
\[
(\ell,q)=1\text{ for any } q|(p-1),\qquad k=\frac{e}{\ell}.
\]
Note that then $f|k$ and also there exists some integer $N$ such that $k|(f^N).$
Since $\gamma$ is a generator for the group $\mathbb{F}_{p}^{\times}/\mathbb{F}_{p}^{\times f}$, then for any $y\in\mathbb{F}_{p}^{\times}$ there exists an unique $i\in\{0,\ldots,f-1\}$ and $z_1\in\mathbb{F}_{p}^{\times}$ such that $y=\gamma^{i}z_{1}^{f}$. On the other hand, we can apply the same principle to $z_{1}$,  finding an unique $j\in\{0,\ldots,f-1\}$ and $z_{2}\in\mathbb{F}_{p}^{\times}$ such that $z_1=\gamma^{j}z_{2}^{f}$. Thus, $y=\gamma^{i}z_{1}^{f}=\gamma^{i}(\gamma^{j}z_{2}^{f})^{f}=\gamma^{i}z_{2}^{f^{2}}$. Iterating this process $N$ times, we can find $z_{N}\in\mathbb{F}_{p}^{\times}$ such that $y=\gamma^{i}z_{N}^{f^{N}}$ with $k|f^{N}$. Then,
$y=\gamma^{i}(z_{N}^{f^{N}/k})^{k}$. On the other hand, since $(\ell,p-1)=1$, we have that $z_{N}^{f^{N}/k}=z^{\ell}$ for some $z\in\mathbb{F}_{p}^{\times}$, so that $y=\ga^i z^{\ell k} = \ga^i z^e$ and this proves the claim.  Moreover, note that once we have obtained $z$ such that $y=\ga^iz^e$ then we can multiply $z$ by any $f$-th root of unity, so that there are $f$ such values $z$.

Next, for any $i\in\{0,\ldots,f-1\}$ we can consider the map
\[
\varphi_{i}:\begin{matrix}
W_{i}& \longrightarrow & W \\ (z,\bfa) &\mapsto & (\gamma^{i}z^{e},\bfa).
\end{matrix}
\]
From this, we deduce that if $(y,\bfa)$ is in the image of $\varphi_{i}$ then
\[
|\varphi_{i}^{-1}(y,\bfa)|=
\begin{cases}
f               &\text{if $y\neq 0$}\\
1               &\text{if $y=0$}.
\end{cases}
\]
On the other hand, if $(0,\bfa)\in W$, then $(0,\bfa)\in W_{i}$ for each of $i=0,\ldots,f-1$. Then the result follows.
\end{proof}

When we apply Lemma \ref{lemma_split_by_root} it will be convenient to treat all cases analogously as $i$ varies; to do so we will employ the following lemma.
\begin{lemma}\label{lemma_isom}
Fix $e \geq 1$ and recall $F(Y,\Xbf)$ from (\ref{F_dfn}).
Let $p$ be a prime, and let $\bfu\in\overline{\mathbb{F}}_{p}^{n}$. Then  for any $\alpha\in\overline{\F}_{p}^{\times}$ the variety
$V(F(\alpha Z^{e},\X))\cap V(\langle\X,\bfu\rangle)\subset\mathbb{P}_{\overline{\mathbb{F}}_{p}}^{n}$
is isomorphic to $V(F(Z^{e},\X))\cap V(\langle\X,\bfu\rangle)\subset\mathbb{P}_{\overline{\mathbb{F}}_{p}}^{n}$.  In particular, for $\ubf = \mathbf{0},$ we conclude $V(F(\alpha Z^{e},\X)) \subset\mathbb{P}_{\overline{\mathbb{F}}_{p}}^{n}$
is isomorphic to $V(F(Z^{e},\X)) \subset\mathbb{P}_{\overline{\mathbb{F}}_{p}}^{n}$.
\label{lem : iso}
\end{lemma}
\begin{proof}
Let $\beta\in\overline{\mathbb{F}}_{p}^{\times}$ be such that $\beta^{e}=\alpha$. Then the change of variables $(Z,\X)\mapsto (\beta Z,\X)$ induces an isomorphism between $V(F(Z^{e},\X))\cap V(\langle\X,\bfu\rangle)$ and $V(F(\alpha Z^{e},\X))\cap V(\langle\X,\bfu\rangle)$.

\end{proof}

\subsection{Proof of Proposition \ref{prop_g_sum}}
We are now ready to prove our main result of this section, Proposition \ref{prop_g_sum}.
In the following, we denote $f=(e,p-1)$. An application of Lemma $\ref{lem : homogen}$ leads to
\begin{equation}
g(\bfu,p)=-\delta_{\bfu=\mathbf{0}}p^{n}+\frac{1}{f}\sum_{i=0}^{f-1}\sum_{\substack{(z,\bfa ) \in W_{i}}}e_{p}(\langle\bfa,\bfu\rangle).
\label{eq : sum}    
\end{equation}

\subsubsection{Type zero case} Assume $\bfu \con \boldsymbol{0} \modd{p}$. The right hand side of $(\ref{eq : sum})$ becomes
\[
g(\mathbf{0},p)=- p^{n}+\frac{1}{f}\sum_{i=0}^{f-1}\sum_{\substack{(z,\bfa ) \in W_{i}}}1=-p^{n}+\frac{1}{f}\sum_{i=0}^{f-1}|W_{i}|.
\]
By definition, for any $i=0,\ldots,f-1$ the set $W_{i}$ is the set of the $\mathbb{F}_{p}$-points on the affine variety $V(F(\gamma^{i}Z^{e},\X))\subset\mathbb{A}^{n+1}_{\mathbb{F}_{p}}$. By hypothesis,    $p$ is of good reduction for $V(F(Z^{e},\X))$,  so  $V(F(Z^{e},\X))\subset\mathbb{P}^{n}_{\overline{\mathbb{F}}_{p}}$ is nonsingular. Then by Lemma \ref{lem : iso}, we have that $V(F(\gamma^{i}Z^{e},\X))\subset\mathbb{P}^{n}_{\overline{\mathbb{F}}_{p}}$ is a nonsingular variety for each $i=0,\ldots,f-1$ (and in particular is absolutely irreducible over $\overline{\F}_p$), and certainly $V(F(\gamma^{i}Z^{e},\X))$ is defined over $\F_p.$
Thus the Lang-Weil bound \cite{LanWei54} implies that  (counting projectively)
\[
 |V(F(\gamma^{i}Z^{e},\X))(\mathbb{F}_{p})|=p^{n-1}+O_{m,e,d}(p^{n-1-1/2})\qquad\text{for each }i=0,\ldots,f-1,
\]
so that 
$|W_{i}|=p^{n}+O_{m,e,d,}(p^{n-1/2})$ for each $i=0,\ldots,f-1.$
Thus we may conclude that $g(\boldsymbol{0},p)\ll p^{n-1/2}$.

\subsubsection{Good/Bad case} Assume $\bfu\neq \boldsymbol{0} \modd{p}$; we may initially argue the good and the bad cases together. The right hand side of $(\ref{eq : sum})$ becomes
\[
g(\bfu,p)=\frac{1}{f}\sum_{i=0}^{f-1}\sum_{\substack{(z,\bfa ) \in W_{i}}}e_{p}(\langle\bfa,\bfu\rangle).
\]
In either the good or the bad case, it suffices  to estimate each sum
\[
g_{i}(\bfu,p)=\sum_{\substack{(z,\bfa ) \in W_{i}}}e_{p}(\langle\bfa,\bfu\rangle), \qquad \text{for $i=0,..,f-1$}.
\]
  First we prove that for any $\alpha\in\mathbb{F}_{p}^{\times}$,  $g_{i}(\bfu,p)=g_{i}(\alpha \bfu,p)$. Indeed  
\[
\begin{split}
g_{i}(\alpha  \bfu,p)&=\sum_{\substack{(z,\bfa ) \in W_{i}}}e_{p}(\langle\bfa,\alpha  \bfu\rangle) =\sum_{\substack{(z,\bfa ) \in \mathbb{F}_{p}^{n+1}\\F(\gamma^{i}z^{e},\bfa )=0}}e_{p}(\langle\bfa,\alpha \bfu\rangle)\\&= \sum_{\substack{(z,\bfa ) \in \mathbb{F}_{p}^{n+1}\\F(\gamma^{i}z^{e},\bfa )=0}}e_{p}(\langle\alpha  \bfa,\bfu\rangle) = \sum_{\substack{(t,\bfb ) \in \mathbb{F}_{p}^{n+1}\\\overline{\alpha}^{med}F(\gamma^{i}t^{e},\bfb )=0}}e_{p}(\langle \bfb,\bfu\rangle)\\&= \sum_{\substack{(t,\bfb ) \in \mathbb{F}_{p}^{n+1}\\F(\gamma^{i}t^{e},\bfb )=0}}e_{p}(\langle \bfb,\bfu\rangle) = g_{i}(\bfu,p),
\end{split}
\]
where in the fourth step we use the change of variables $(z,\bfa)=(\overline{\alpha}t,\overline{\alpha} \bfb)$, for $\al \overline{\al} \con 1 \modd{p}$. Hence,
\[
\begin{split}
(p-1)g_{i}(\bfu,p)&=\sum_{\alpha\in\mathbb{F}_{p}^{\times}}g_{i}(\alpha \bfu,p)\\& =\sum_{\alpha\in\mathbb{F}_{p}^{\times}}\sum_{\substack{(z,\bfa ) \in \mathbb{F}_{p}^{n+1}\\F(\gamma^{i}z^{e},\bfa )=0}}e_{p}(\langle\bfa,\alpha \bfu\rangle)\\
&=\sum_{\substack{(z,\bfa ) \in \mathbb{F}_{p}^{n+1}\\F(\gamma^{i}z^{e},\bfa )=0}}\sum_{\alpha\in\mathbb{F}_{p}^{\times}}e_{p}(\alpha \langle\bfa,\bfu\rangle)
    =\sum_{\substack{(z,\bfa ) \in \mathbb{F}_{p}^{n+1}\\F(\gamma^{i}z^{e},\bfa )=0}}\sum_{\alpha\in\mathbb{F}_{p}}e_{p}(\alpha \langle\bfa,\bfu\rangle) - \sum_{\substack{(z,\bfa ) \in \mathbb{F}_{p}^{n+1}\\F(\gamma^{i}z^{e},\bfa )=0}} 1\\
&=p(p-1)\cdot |(V(F(\gamma^{i}Z^{e},\X))\cap V(\langle\bfu,\X\rangle))(\mathbb{F}_{p})|-(p-1)\cdot |V(F(\gamma^{i}Z^{e},\X)(\mathbb{F}_{p})|+(p-1),
\end{split}
\]
where in the last step we have passed to counting points over $\F_p$ in the projective sense.
Applying \cite[Appendix by N. Katz, Theorem $1$]{Hoo91}, we have that
\[
\begin{split}
&|V(F(\gamma^{i}Z^{e},\X))(\mathbb{F}_{p})|=\sum_{j=0}^{n-1}p^{j}+O_{n,m,e,d}(p^{\frac{n+\delta_{i}}{2}})\\&|(V(F(\gamma^{i}Z^{e},\X))\cap V(\langle\bfu,\X\rangle))(\mathbb{F}_{p})|=\sum_{j=0}^{n-2}p^{j}+O_{n,m,e,d}(p^{\frac{n-1+\delta_{i,\bfu}}{2}}),
\end{split}
\]
where $\delta_{i}=\dim (\sing (V(F(\gamma^{i}Z^{e},\X))$ and $\delta_{i,\bfu}=\dim(\sing (V(F(\gamma^{i}Z^{e},\X))\cap V(\langle\bfu,\X\rangle)))$.

On the other hand, Lemma $\ref{lem : iso}$ implies that $\delta_{i}=\delta_{0}$ and $\delta_{i,\bfu}=\delta_{0,\bfu}$ for each $i$. Moreover, $\delta_{0}=-1$ since we are assuming that $p$ is of good reduction for $V(F(Z^{e},\X))$. Thus, we obtain
\beq\label{Katz_app}
g_{i}(\bfu,p)=O(p^{\frac{n+1+\delta_{0,\bfu}}{2}}),
\eeq
 with an implicit constant depending only on $n,m,e,d$.
Finally, by (\ref{sing_dim}),   
\[
\delta_{0,\bfu}=
\begin{cases}
0          &\text{if $V(\langle\bfu,\X\rangle)$ is tangent to $V(F(Z^{e},\X))$}\\
-1         &\text{otherwise},
\end{cases}
\]
and this completes the proof of the good and bad cases in Proposition \ref{prop_g_sum}.

\begin{remark}\label{remark_sing_dim}
This remark justifies (\ref{sing_dim}).   Let $V=V(H(\X)) \subset \mathbb{P}^n$ be a nonsingular hypersurface and $L=V(\langle \bfa,   \X\rangle)$ be a hyperplane. We may suppose without loss of generality that $a_1 \neq 0.$ By the Jacobian criterion,
$\sing (V \intersect L)$ is the set of points on the intersection  $V \intersect L$ for which the $ (n+1)\times 2$ matrix with columns $\nabla H$ and $\abf$ has rank 1. Consequently,  $\sing (V \intersect L)\subset W$ where 
\[ W = V \intersect V(g_2) \intersect \cdots \intersect V(g_n),\]
in which for each $i=2,\ldots, n,$ 
\[
g_i(\Xbf)  = a_1 \frac{\partial H}{\partial X_i}(\X)  - a_i \frac{\partial H}{\partial X_1}(\X).
\]
On the other hand, $W \intersect V(\partial H/\partial X_1) = \sing (V) =\emptyset$ under the hypothesis that $V$ is nonsingular. Consequently, $\dim W \leq 0,$ implying $\dim(\sing(V \intersect L )) \leq 0,$ as desired.
\end{remark}

 \begin{remark}\label{remark_F_hom}
 It is worth remarking what we have gained from the arguments in this section. Briefly, suppose $\bfu \not\con 0 \modd{p}$ and consider
\[
g(\ubf,p) =\sum_{\substack{(y,\bfa ) \in\mathbb{F}_{p}^{n+1}\\ F(y,\bfa)=0}}e_{p}(\langle\bfa,\bfu\rangle).
\]
To work directly with this sum rather than passing through the dissection into the components $W_i$ as we did above, we would first need to homogenize the polynomial $F(Y,\bfx)$, say defining a homogeneous polynomial \[\tilde{F}(T,Y,\Xbf) =
T^{md(e-1)}Y^{md} + \cdots +
T^{m(e-1)}Y^mf_{d-1}(\Xbf) + f_d(\Xbf).
\]
(Here we suppose that $e\geq 2$ for this example.)
Then observe that $[1:0:\ldots:0]$ is a singular point on $V(\tilde{F}(T,Y,\Xbf)) \subset \mathbb{P}^{n+1}.$
Consequently, if one proceeded to estimate $g(\ubf,p)$, roughly analogous to the approach in (\ref{Katz_app}), by counting points on the complete intersection described by $V(\tilde{F}(T,Y,\Xbf))\intersect V(\langle \bfu,\Xbf\rangle) \intersect V(T=1)$, the role of $\del_{0,\ubf}$ in the exponent is now played by a dimension that is always at least $0$, ultimately leading to a result that is larger by a factor of $p^{1/2}$ than the results we obtain in Proposition \ref{prop_g_sum}.
 \end{remark}

\subsection{Choice of the sieving set}\label{sec_sieve_set}
We can now continue the discussion initiated in \S \ref{sec_sieve_first}, and choose the sieving set. 
We suppose that $Q=B^\kappa$ for some $1/2 \leq \kappa \leq 1$ to be chosen later (see (\ref{Q_choice})). 
We choose the sieving set 
\[ \Pcal \subset  [Q,2Q] \]
 comprised of all primes  in this range such that 
(i) $p \con 1 \modd{m}$ (recalling $m \geq 2$), and (iii') the reduction $V(F(Z^e,\X)) \subset \mathbb{P}_{\overline{\F}_p}^n$ is nonsingular. 

 By the Siegel--Walfisz theorem on primes in arithmetic progressions,  there are $\gg_{m} Q/\log Q$ primes such that $p\con 1 \modd{m}$  in any dyadic region $[Q,2Q],$ for all $Q \gg_m 1$ sufficiently large, which we assume is a condition met henceforward. We recall from (\ref{few_primes}) that at most $O_{m,e,d}(\log \|F\|)$ primes fail (iii'). We henceforward assume that 
\beq\label{Q_lower}
Q \gg_{m,e,d} (\log \|F\|)(\log \log \|F\|)
\eeq
for an appropriately large implied constant, so that consequently
\beq\label{P_lower}
P = |\Pcal| \gg_{m} Q/\log Q  - C_{m,e,d}(\log \|F\|) \gg_{m,e,d}Q/\log Q.
\eeq
 When we finally choose $Q$ as a power of $B$, (\ref{Q_lower}) will impose a lower bound on $B$; we defer this to (\ref{Q_choice}).

\section{Estimating the main sieve term: the bad-bad case}\label{sec_badbad}
This section is the technical heart of the paper. We show how to bound the most difficult contribution to the sieve, which occurs when $\ubf$ is bad with respect to two primes $p \neq q \in \Pcal$. (We reserve the treatment of all other cases, when $\ubf$ is either type zero, or good with respect to at least one of these primes, to \S \ref{sec_concluding}; these remaining cases are significantly easier.)

We recall from the sieve lemma, Lemma \ref{lemma_sieve}, that $\Scal(F,B)$ is bounded above by a sum of three terms. The first two terms can be bounded simply:
\beq\label{S_simple}
\sum_{\bfk:f_{d}(\bfk)=0}W(\bfk)+ \frac{1}{P}\sum_{\bfk}W(\bfk) \ll  B^{n-1}+B^nP^{-1}.
\eeq
 Here the first term follows from the Schwartz-Zippel trivial bound $\ll_{n,e,d}B^{n-1}$ for the number of zeroes of $f_d$ with $\kbf \in \supp(W)$, since $f_d \not\con 0$ (see e.g. \cite[Theorem 1]{HB02}, which as mentioned before has a method of proof that applies even if $f_d$ is not absolutely irreducible).
 We will call the remaining, third, term on the right-hand side of the sieve lemma the main sieve term.

Now we are ready to estimate the main sieve term, which after an application of Poisson summation inside the definition (\ref{T_dfn}) of $T(p,q)$ is
\begin{align}
\frac{1}{P^{2}}\sum_{\substack{p,q\in\mathcal{P}\\p\neq q}}|T(p,q)|&=\frac{1}{P^{2}}\sum_{\substack{p,q\in\mathcal{P}\\p\neq q}}\left(\frac{1}{pq}\right)^{n}\left|\sum_{\bfu}\hat{W}\left(\frac{\bfu}{pq}\right)  g(\bfu,pq) \right| \nonumber \\
&
\ll\frac{1}{P^{2}Q^{2n}}\sum_{\substack{p,q\in\mathcal{P}\\p\neq q}}\sum_{\bfu}\left|\hat{W}\left(\frac{\bfu}{pq}\right) g(\bfu,pq)\right|. \label{sieve_Tpq}
\end{align}
We will apply Proposition \ref{prop_g_sum} to bound $g(\ubf,pq)$, according to the ``type'' of $\ubf$ modulo $p$ and $q$, respectively; this leads to cases we can abbreviate as zero-zero, zero-good, zero-bad, good-good, good-bad, and bad-bad. Unsurprisingly, the greatest difficulty is to bound the contribution of the bad-bad case, and we focus on this first, returning to the other cases in \S \ref{sec_concluding}.

Recall that $W$ is a   non-negative function with $W(\bfu) = w(\bfu/B)$ for an infinitely differentiable, non-negative function $w$ that is $\con 1$ on $[-1,1]$ and vanishes outside of $[-2,2]$. Thus $\hat{W}(\bfu) = B^n \hat{w} (B \bfu)$ and $\hat{w}(\bfu)$ has rapid decay in $\bfu$, so that
\beq\label{WB_bound}
|\hat{W}(\bfu)|\ll B^n \prod_{i=1}^{n}\left(1+|u_{i}|B\right)^{-M}
\eeq
for any $M \geq 1$; we will for example specify a lower bound on $M$ at (\ref{M_choose}) and can certainly always assume $M \geq 2n$.
 In particular, we will later apply the fact that for any $B,L \geq 1$,
\beq\label{WB_bound_sum} \sum_{\ubf \in \Z^n} |\hat{W}(\ubf/L)| \ll \max\{B^n,L^n\}.
\eeq

\subsection{The dual variety}\label{sec_dual}
To consider any bad case, it is useful to consider certain facts about the dual variety. 
Recall that $m \geq 2$ and   $d, e \geq 1$, and
\beq\label{F_dfn_repeat}
F(Y,\X)=Y^{md}+Y^{m(d-1)}f_{1}(\X)+\ldots+f_{d}(\X),
\eeq
in which for each $1 \leq i \leq d$, $f_i$ is a polynomial in $\Z[X_1,\ldots, X_n]$ with  $\deg f_{i}=m\cdot e\cdot i$.
By hypothesis, the variety defined by $F(Y,\X)=0$ in weighted projective space, denoted  $V(F(Y,\Xbf)) \subset \mathbb{P}_{\mathbb{C}}(e,1,\ldots,1),$ is nonsingular.
 Recall from \S \ref{sec_smooth} that $V(F(Y,\Xbf)) \subset \mathbb{P}_{\mathbb{C}}(e,1,\ldots,1)$ is nonsingular if and only if $V(F(Z^e,\Xbf)) \subset \mathbb{P}_{\mathbb{C}}^n$ is nonsingular.
The dual variety $V^{*}=V(F(Z^{e},\X))^{*}\subset\mathbb{P}^{n}_\mathbb{C}$  of a hypersurface is a hypersurface. 
We denote by
\beq\label{G_dfn}
G(U_{Y},U_{1},\ldots,U_{n})  
\eeq
  the irreducible homogeneous  polynomial such that $V(G)=V^{*}$ (see e.g. \cite[Prop. 11.2, Appendix]{BCLP23}).
Recall that $\deg F(Z^{e},\X) = mde $; by \cite[Prop. 2.9]{EisHar16},
\[\deg G = mde (mde-1)^{n-1} \geq 2.\]
 In our analysis of the bad-bad case in   \S \ref{sec_bad_bad}, our strategy is to divide our analysis depending on whether $\ubf$ has the property $G(0,\ubf)\neq0$ or $G(0,\ubf)=0$. In the first case, we now show via an explicit constructive argument that
\beq\label{prime_log_F}
|\{p : \text{$\mathbf{u}$ is bad modulo $p$}\}| \ll_{n,m,e,d} \log (\|F\| \|\ubf\|).\eeq

Let us prove this.
A given $\bfu$ has the property $G(0,\ubf) \neq 0$ if and only if the hyperplane $V(\langle\bfu, \X\rangle )\subset\mathbb{P}_\mathbb{C}^{n}$ is not tangent to $V(F(Z^{e},\X))\subset\mathbb{P}_\mathbb{C}^{n}$; that is, if and only if     for any $[z:\x]\in V(F(Z^{e},\X))\cap V(\langle\X,\bfu\rangle )$, the matrix
\beq\label{tan_matrix}
\begin{pmatrix}
\frac{\partial F}{\partial Z}(z^{e},\x) & 0 \\
\frac{\partial F}{\partial X_{1}}(z^{e},\x)  & u_{1}\\
\vdots\\
\frac{\partial F}{\partial X_{n}}(z^{e},\x)  & u_{n}
\end{pmatrix}
\eeq
has maximal rank (i.e. at least one $2 \times 2$ minor is nonvanishing).  
Now define $n+2$ polynomials in $Z,X_1,\ldots,X_n$, with integral coefficients (depending on $\bfu$) as follows: set 
\[H_{0,\bfu}(Z,\X)=F(Z^{e},\X), \qquad H_{n+1,\bfu}(Z,\X)=\langle\X,\bfu\rangle,
\]
and
for $1 \leq i \leq n$ set
\[
H_{i,\bfu}(Z,\X)=
\begin{cases}
\det\begin{pmatrix}
\frac{\partial F}{\partial Z}(z^{e},\x) & 0\\
\frac{\partial F}{\partial X_{1}}(z^{e},\x) & u_{1}\end{pmatrix}&\text{for $i= 1$}\\
\det\begin{pmatrix}\frac{\partial F}{\partial X_{i-1}}(z^{e},\x) & u_{i-1}\\
\frac{\partial F}{\partial X_{i}}(z^{e},\x) & u_{i}
\end{pmatrix}		&\text{for $2\leq i\leq n.$}
\end{cases}
\] 
Then define the resultant (see \cite[Ch. 13]{GKZ08})
\beq\label{R_dfn}
R (\bfu)=\text{Res}(H_{0,\bfu},H_{1,\bfu},\ldots,H_{n+1,\bfu}). 
\eeq
The following are all equivalent:
\begin{enumerate}
    \item $\bfu$ has the property that $V(\langle \bfu,\X\rangle)$ is tangent to $V(F(Z^e,\X))$
    \item for some $[z:\x]\in V(F(Z^{e},\X))\cap V(\langle\X,\bfu\rangle )$, (\ref{tan_matrix}) has rank $<2$
    \item the polynomials $H_{i,\bfu}(Z,\X)$   (for $0  \leq i \leq n+1$) share a common (nonzero) root
    \item $R(\bfu) = 0$.
\end{enumerate}  

Now we consider the analogues of these statements for each $p$. Fix a prime $p$. For a polynomial $L\in\mathbb{Z}[\bfU]$, let $\overline{L}$ denote its reduction modulo $p$.  By definition, $\bfu$ is bad modulo $p$ precisely when $\overline{H}_{i,\bfu}$ (for $0 \leq i \leq n+1$) have a common nontrivial root modulo $p$, that is if and only if $p | \text{Res}(\overline{H}_{0,\bfu},\ldots,\overline{H}_{n+1,\bfu})$. By \cite[Section IV]{Cha93}, as a polynomial in $\bfU,$
\[
\text{Res}(\overline{H}_{0,\bfU},\ldots,\overline{H}_{n+1,\bfU})= \overline{R}(\bfU),
\]
where $R$ is defined as in (\ref{R_dfn}). (That is, the resultant of the reductions modulo $p$ is the reduction modulo $p$ of the resultant.)
Thus for each $\bfu$ such that $G(0,\bfu) \neq 0$ so that $R(\bfu) \neq 0,$ we can conclude that 
\[ |\{ p : \text{$ \bfu $ is bad modulo $p$}\}| = \omega ( \mathrm{Res}(H_{0,\bfu}, \ldots, H_{n+1,\bfu})),\]
 where $\omega(r)$ indicates the number of distinct prime divisors of an integer $r$; we recall in particular that $\om(r) \ll \frac{\log r}{\log \log r}$.
 By \cite{GKZ08}[Prop. 1.1, Ch. 13], the resultant is a homogeneous polynomial in the coefficients of the forms $H_{0,\bfu}, \ldots, H_{n+1,\bfu}$ (with degree bounded in terms of $n,m,e,d$). Thus, for every value of $\ubf$ such that $G(0,\ubf) \neq 0$ so that $\mathrm{Res}(H_{0,\bfu}, \ldots, H_{n+1,\bfu})$ is a nonzero integer,
\beq\label{omega_Res}
  \omega ( \mathrm{Res}(H_{0,\bfu}, \ldots, H_{n+1,\bfu})) \ll_{n,m,e,d} \log (\|F\|  \| \bfu \|).
  \eeq
  
 Finally,   if $G(0,\bfu)=0$, then the hyperplane $V(\langle \ubf,\X\rangle) \subset \mathbb{P}_\C^n$ is tangent to $V(F(Z^e,\X)) \subset \mathbb{P}_\C^n$ so that (\ref{tan_matrix}) has rank 1 over $\C$; consequently $\bfu$ is bad for all primes $p$. Thus in this latter case, we will instead focus on showing there are sufficiently few solutions to $G(0,\bfu)=0$.

\begin{remark} It is a common occurrence that one requires the fact that there are ``quite few'' primes of bad reduction for a variety of the form $\mathcal{V} \intersect \{u_0X_0 + \cdots u_nX_n=0\}$ for some variety $\mathcal{V}$ and parameter $(u_0,u_1,\ldots,u_n)$ of interest, in this case $V(G)$ with $G$ describing the dual of $F$, and $u_0=0$. The fact that our result (\ref{prime_log_F}) depends only logarithmically on $\|F\|$   is important for our ultimate deduction that the implicit constant in Theorem \ref{thm_main} is independent of $\|F\|$; see the application in \S \ref{sec_G_nonzero}. This motivated the explicit argument we gave above. Alternatively, we thank Per Salberger for pointing out that the useful references \cite[pp. 95-98]{CLO05} and \cite{Dem12}  also provide  similar constructions leading to explicit results of the form (\ref{omega_Res}) and hence (\ref{prime_log_F}). We remark that if we did not require logarithmic dependence on $\|F\|$, one could apply a result such as \cite[Prop. 11.5(3), Appendix] {BCLP23} to conclude immediately that for all sufficiently large primes (in an inexplicit sense), $\ubf$ is bad modulo $p$ precisely when $p|G(0,\ubf)$ (so that $|\{ p: \text{$\ubf$ is bad modulo p}\}|\ll_G \log \|\ubf\|$ when $G(0,\bfu) \neq 0$), but with  dependence on $G$ and hence on $F$ that has not been made explicit, and so does not immediately suffice for our application.
\end{remark}

\subsection{Bad-bad case}\label{sec_bad_bad}
We use the above facts to control the contribution of the bad-bad case to the sieve, which by Proposition \ref{prop_g_sum} is bounded by
\begin{equation}
\begin{split}
\frac{1}{P^2 Q^{2n}}\sum_{\substack{p,q\in\mathcal{P}\\p\neq q}}\sum_{\substack{\bfu\in\mathbb{Z}^{n}\\\bfu\text{ bad mod }p\\\bfu\text{ bad mod }q}}\left|\hat{W}\left(\frac{\bfu}{pq}\right) g(\bfu,pq)\right|&\ll \frac{Q^{n+1}}{P^{2}Q^{2n}}\sum_{\substack{p,q\in\mathcal{P}\\p\neq q}}\sum_{\substack{\bfu\in\mathbb{Z}^{n}\\\bfu\text{ bad mod }p\\\bfu\text{ bad mod }q}}\left|\hat{W}\left(\frac{\bfu}{pq}\right)\right|.
\end{split}
\label{badbad}
\end{equation}
We start by exchanging the order of summation between $\bfu$ and the primes $p,q$, and then splitting the sum as
\[
\sum_{\substack{\bfu\in\mathbb{Z}^{n}}}\sum_{\substack{p,q\in\mathcal{P}\\p\neq q\\\bfu\text{ bad mod }p\\\bfu\text{ bad mod }q}}\left|\hat{W}\left(\frac{\bfu}{pq}\right)\right|\\
    =\sum_{\substack{\bfu\in\mathbb{Z}^{n}\\G(0,\bfu)=0}}\sum_{\substack{p,q\in\mathcal{P}\\p\neq q\\\bfu\text{ bad mod }p\\\bfu\text{ bad mod }q}} +\sum_{\substack{\bfu\in\mathbb{Z}^{n}\\G(0,\bfu)\neq 0}}\sum_{\substack{p,q\in\mathcal{P}\\p\neq q\\\bfu\text{ bad mod }p\\\bfu\text{ bad mod }q}} .
\]

In this section, we will prove that the contribution from $G(0,\bfu) \neq 0$ is
\beq\label{badbad_unonzero}
\sum_{\substack{\bfu\in\mathbb{Z}^{n}\\G(0,\bfu)\neq 0}}\sum_{\substack{p,q\in\mathcal{P}\\p\neq q\\\bfu\text{ bad mod }p\\\bfu\text{ bad mod }q}}\left|\hat{W}\left(\frac{\bfu}{pq}\right)\right| \ll_{n,m,e,d}Q^{2n}(\log B)^2.
\eeq
On the other hand, we will prove that the contribution from $G(0,\bfu) = 0$ is
\beq\label{Wpq_bad}
 \sum_{\substack{\bfu\in\mathbb{Z}^{n}\\G(0,\bfu)=0}}\sum_{\substack{p,q\in\mathcal{P}\\p\neq q\\\bfu\text{ bad mod }p\\\bfu\text{ bad mod }q}}\left|\hat{W}\left(\frac{\bfu}{pq}\right)\right|\ll_{\varepsilon} P^2\left(Q^{2n}  B^{-\al(M-1)}+ B^n \left(\frac{Q^{2}}{B^{1-\alpha}}\right)^{n-2+\frac{1}{3}+\varepsilon}\right),\eeq
 for a small $0<\al<1$ of our choice, and any $\ep>0$. 
 Once we have proved these two inequalities, we will wrap up the contribution of the bad-bad case in \S \ref{sec_badbad_wrapup}.
 
\subsubsection{The case $G(0,\bfu) \neq0$}\label{sec_G_nonzero}
Proving (\ref{badbad_unonzero}) is quite simple; by the decay (\ref{WB_bound}) for $\hat{W}$ and the bound (\ref{omega_Res}) for counting $p,q,$
\begin{align*}
\sum_{\substack{\bfu\in\mathbb{Z}^{n}\\G(0,\bfu)\neq 0}}\sum_{\substack{p,q\in\mathcal{P}\\p\neq q\\\bfu\text{ bad mod }p\\\bfu\text{ bad mod }q}}\left|\hat{W}\left(\frac{\bfu}{pq}\right)\right|
&\ll B^n \sum_{\substack{\bfu\in\mathbb{Z}^{n}\\G(0,\bfu)\neq 0}}\prod_{i=1}^{n}\left(1+\frac{B|u_{i}|}{Q^{2}}\right)^{-M}\omega(R(\bfu))^{2} \nonumber\\
&\ll  B^n \sum_{\bfu\in\mathbb{Z}^{n}}\prod_{i=1}^{n}\left(1+\frac{B|u_{i}|}{Q^{2}}\right)^{-M}(\log ( \|F\| \|\ubf\|))^{2} \nonumber\\
&\ll_{n,m,e,d}  Q^{2n}  (\log B)^{2}. 
\end{align*}
Here we have used  the fact that $Q=B^\kappa$ with $1/2 \leq \kappa \leq 1$ (so that $Q^{2n} \gg B^n$), and the  fact from Lemma \ref{lemma_F_small} that in the only case we need to consider, $\log \|F\| \ll_{m,e,d}  \log B.$ This proves (\ref{badbad_unonzero}) with an implied constant independent of $\|F\|$.

\subsubsection{The case $G(0,\bfu)=0$} 
Proving (\ref{Wpq_bad}) is a key novel aspect of our proof. Note that if $G(0,\bfu)=0$, then $\bfu$ is bad mod $p$ for all $p\in \mathcal{P}$. Then
\beq\label{badbad_vanish}
\sum_{\substack{\bfu\in\mathbb{Z}^{n}\\G(0,\bfu)=0}}\sum_{\substack{p,q\in\mathcal{P}\\p\neq q\\\bfu\text{ bad mod }p\\\bfu\text{ bad mod }q}}\left|\hat{W}\left(\frac{\bfu}{pq}\right)\right|
	\ll  B^n P^{2} \sum_{\substack{\bfu\in\mathbb{Z}^{n}\\G(0,\bfu)=0}}\prod_{i=1}^{n}\left(1+\frac{B|u_{i}|}{Q^{2}}\right)^{-M}.
\eeq
Let $0<\alpha < 1$ be a parameter to be chosen later and consider the cube
\[\mathcal{C}_{\alpha}=[-Q^{2}/B^{1-\alpha},Q^{2}/B^{1-\alpha}]^{n} \subset \R^n.\]
This is slightly larger than the ``essential support'' of the sum over $\bfu$, so that outside this box we can exploit decay more efficiently.
We will ultimately prove that
\beq\label{W_G_above}
\sum_{\substack{\bfu\in\mathbb{Z}^{n}\\G(0,\bfu)=0}}\prod_{i=1}^{n}\left(1+\frac{B|u_{i}|}{Q^{2}}\right)^{-M} \ll_\ep Q^{2n}B^{-n}B^{-\al(M-1)} +  \left( \frac{Q^2}{B^{1-\al}}\right)^{n-2+1/3+\ep},
\eeq
for any $\ep>0.$
We split the sum as  
\beq\label{two_remaining}
\sum_{\substack{\bfu\in\mathcal{C}_{\alpha} \intersect \Z^n\\G(0,\bfu)=0}}\prod_{i=1}^{n}\left(1+\frac{B|u_{i}|}{Q^{2}}\right)^{-M}+\sum_{\substack{\bfu\notin\mathcal{C}_{\alpha} \intersect \Z^n\\G(0,\bfu)=0}}\prod_{i=1}^{n}\left(1+\frac{B|u_{i}|}{Q^{2}}\right)^{-M}.
\eeq
In the second sum in (\ref{two_remaining}), we can exploit decay:
\[
\sum_{\substack{\bfu\notin\mathcal{C}_{\alpha}\\G(0,\bfu)=0}}\prod_{i=1}^{n}\left(1+\frac{B|u_{i}|}{Q^{2}}\right)^{-M}
\ll\sum_{j=1}^{n}\sum_{\substack{\bfu\in\mathbb{Z}^{n}\\G(0,\bfu)=0\\|u_{j}|>Q^{2}/B^{1-\alpha}}}\prod_{i=1}^{n}\left(1+\frac{B|u_{i}|}{Q^{2}}\right)^{-M} 
\ll \left(\frac{Q^{2}}{B}\right)^n  \frac{1}{B^{\al (M-1)}} .
\]
The contribution of these $\bfu$ to (\ref{badbad_vanish}) is thus 
$
\ll Q^{2n} P^2 B^{-\al(M-1)}
$
 for $0<\al<1$ and any $M \geq 2n$; this contributes the first term in (\ref{Wpq_bad}).

It remains to deal with the first sum appearing on the right hand side of (\ref{two_remaining}), summing over $\bfu \in \mathcal{C}_\al$ such that $G(0,\bfu)=0$. Here we show that there are few solutions to $G(0,\bfu)=0$.
Recall  the definition of the form $G$ from \S \ref{sec_dual}.
 Consider $V(G(0,\bfU)) \subset \mathbb{P}^{n-1}$ defined by $G(0,\mathbf{U})=0$ as a function of $\Ubf$. (First notice that $G(0,\Ubf)$ is not identically zero; indeed, if it were then we would conclude that $   \{ U_Y =0\} \subset \{G(U_Y,U_1,\ldots,U_n)=0\}$. Recalling that $G(U_Y,\Ubf)$ is irreducible, both these projective varieties have dimension $n-1$ so that in fact we must have $\{G=0\} =\{ U_Y =0\}$. But this is impossible, since $G$ has degree $>1$.) 
Thus $V(G(0,\Ubf))\subset \mathbb{P}_\C^{n-1}$ is a projective variety  of dimension $n-2$  and $\deg G(0,\Ubf) = \deg G(U_Y,\Ubf) \geq 2$.
Moreover, let us decompose $G(0,\bfU)$ into irreducible components, i.e. by writing \beq\label{G_ell}
G(0,\bfU)=\prod_{\ell=1}^{L}G_{\ell}(\bfU),
\eeq
where $G_{\ell}(\bfU)$ is an irreducible polynomial for each $\ell\leq L$ (and $L \ll_{n,m,e,d}1$). Set $d_{\ell}:=\deg G_{\ell}$. We have
\[
\sum_{\substack{\bfu\in\mathcal{C}_{\alpha} \intersect \Z^n\\G(0,\bfu)=0}}\prod_{i=1}^{n}\left(1+\frac{B|u_{i}|}{Q^{2}}\right)^{-M}\leq\sum_{\substack{\bfu\in\mathcal{C}_{\alpha} \intersect \Z^n\\G(0,\bfu)=0}}1\leq\sum_{\ell=1}^{L}\sum_{\substack{\bfu\in\mathcal{C}_{\alpha} \intersect \Z^n\\G_{\ell}(\bfu)=0}}1.
\]
In the next section, we shall prove:
\begin{prop}\label{prop_G_no_linear}
Let $n \geq 3$. For the homogeneous polynomial $G(U_Y,U_1,\ldots,U_n) \in \C[U_Y,U_1,\ldots,U_n]$ defined in (\ref{G_dfn}), $G(0,U_1,\ldots,U_n)$ contains no linear factor, that is, we cannot write $G(0,\bfU)=L(\bfU)\tilde{H}(\bfU)$ for any linear form $L(\bfU) \in \C[U_1,\ldots,U_n].$
\end{prop}

\begin{remark}\label{remark_no_low_rank}
As a consequence of Proposition \ref{prop_G_no_linear},  $G(0,U_1,\ldots,U_n)$ contains no factor in one or two variables. For suppose that in the notation of (\ref{G_ell}) some factor $G_\ell(\bfU)$ (after an appropriate $GL_n(\C)$ change of variables) can be written as a polynomial $g_1(U_1)$ or $g_2(U_1,U_2)$. Then $g_1(U_1)$ is a monomial, hence a product of linear factors, contradicting the proposition. Alternatively,  any form $g_2(U_1,U_2)$     factors over $\C$ into homogeneous linear factors in $U_1,U_2$, as  a consequence of the fundamental theorem of algebra applied to $g_2(1,t) \in \C[t]$, followed by noting $g_2(U_1,U_2) = U_1^{\deg g_2}g_2(1,U_2/U_1).$ This again would contradict the proposition. (Since the statement of Proposition \ref{prop_G_no_linear} is false if $n=2$, see Remark \ref{remark_badbad_n2} for an alternative approach for $n=2$.) 
\end{remark}

The crucial point is that Proposition $\ref{prop_G_no_linear}$ implies that for each  $\ell=1,\ldots,L$ the degree $d_\ell \geq 2$ (and $G_\ell$ depends on at least 3 variables).
By \cite[Theorem 2]{HB02}, and \cite[Theorem A]{Pil95}, we have, for any $\ep>0$,
\beq\label{HBP_bounds}
\sum_{\substack{\bfu\in\mathcal{C}_{\alpha} \intersect \Z^n\\G_{\ell}(\bfu)=0}}1\ll_{\varepsilon}
\begin{cases}
\left(\frac{Q^{2}}{B^{1-\alpha}}\right)^{n-2+\varepsilon}                              &\text{if $d_{\ell} =2$}\\
\left(\frac{Q^{2}}{B^{1-\alpha}}\right)^{n-2+\frac{1}{d_{\ell}}+\varepsilon}                 &\text{if $ d_{\ell} >2$}.\\
\end{cases}
\eeq
 Within these results, the implied constant is independent of $\|F\|$ in each case.
In particular, we may conclude that for each $\ell =1,\ldots, L,$
\[
\sum_{\substack{\bfu\in\mathcal{C}_{\alpha} \intersect \Z^n\\G_{\ell}(0,\bfu)=0}}1\ll_{\varepsilon}\left(\frac{Q^{2}}{B^{1-\alpha}}\right)^{n-2+\frac{1}{3}+\varepsilon}.
\]
Thus the total contribution of these terms to (\ref{badbad_vanish}) is 
\[\ll_\ep B^nP^2\left(\frac{Q^{2}}{B^{1-\alpha}}\right)^{n-2+\frac{1}{3}+\varepsilon}.\]
This contributes the second term in (\ref{Wpq_bad}), and hence (\ref{Wpq_bad})  is proved. 

 \subsubsection{Conclusion of the bad-bad sieve term}\label{sec_badbad_wrapup}
From (\ref{badbad_unonzero}) and (\ref{Wpq_bad}) we conclude that   the total contribution of the bad-bad case (\ref{badbad}) to the sieve is
\begin{multline} \frac{Q^{n+1}}{P^2Q^{2n}} \left( Q^{2n} (\log B)^2 +  Q^{2n} P^2 B^{-\al(M-1)}+ B^nP^2\left(\frac{Q^{2}}{B^{1-\alpha}}\right)^{n-2+\frac{1}{3}+\varepsilon}\right) \\ \ll_{\varepsilon'} Q^n \left(  QP^{-2}(\log B)^2+ QB^{-\al (M-1)}+  \left(\frac{B^{\frac{5}{3}+g(\al)+\varepsilon'}}{Q^{\frac{7}{3}+\varepsilon'}}\right)\right),\label{BQfrac}\end{multline}
where $g(\al)=\al(n-\frac{5}{3}+\varepsilon')$, for any $\ep'>0$. 
To simplify the third term above, henceforward we assume $Q=B^\kappa$ with
\beq\label{kappa_assp} 3/4 \leq \kappa \leq 1.
\eeq
Then the above is
\beq\label{after_Q_lower}
 \ll_{\varepsilon'}   Q^{n}(QP^{-2}(\log B)^2+ QB^{-\al (M-1)}+ B^{-\frac{1}{12}+g(\al)+\varepsilon '}),
\eeq
for any $\ep'>0.$
In the first term on the right-hand side, we observe by (\ref{P_lower}) that $P \gg Q/\log Q$ so that 
\[QP^{-2}(\log B)^2 \ll Q^{-1} (\log B)^4 \ll B^{-3/4}(\log B)^4  .\]
In the second term, we can choose $\alpha=\frac{1}{24}(n-\frac{5}{3}+\varepsilon')^{-1}$ so $g(\al)=1/24,$ and set $M\geq \max\{2n,\al^{-1}+1\}$.
 Regarding the third term, so far this is true for any $\ep'>0$; let us take $\ep'=1/100,$ say.
We conclude that
\beq\label{M_choose}
\frac{Q^{n+1}}{P^{2}Q^{2n}}\sum_{\substack{\bfu\in\mathbb{Z}^{n}\\G(0,\bfu)=0}}\sum_{\substack{p,q\in\mathcal{P}\\p\neq q\\\bfu\text{ bad mod }p\\\bfu\text{ bad mod }q}}\left|\hat{W}\left(\frac{\bfu}{pq}\right)\right|
\ll Q^{n}(B^{-3/4}(\log B)^4 +QB^{-1}+ B^{-\frac{1}{24}+\frac{1}{100}})\ll    Q^{n},
\eeq
since $B \geq Q$. The implied constant is independent of $\|F\|.$ (Here we could even obtain a term that is $o(Q^n)$, but this will not change our main theorem, since the good-good contribution to the sieve is $O(Q^n)$.)
This completes the treatment of the bad-bad contribution to the sieve, except for the proof of Proposition \ref{prop_G_no_linear}, which we provide in the next section. Then in \S \ref{sec_concluding} we show that the contributions of all the other types to the sieve are also dominated by  $\ll Q^{n}$, and then conclude the proof of our main theorem.

 \begin{remark}[The case $n=2$]\label{remark_badbad_n2}
The method of this paper applies for $n=2$ up until Proposition \ref{prop_G_no_linear}; arguing as in Remark \ref{remark_no_low_rank} shows that $G(0,U_1,U_2)$ factors over $\C$ into homogeneous linear factors in $U_1,U_2$, so that proposition is false for $n=2$.
Thus in the nomenclature of (\ref{G_ell}), each degree $d_\ell=1$, and the estimate  
(\ref{HBP_bounds}) is replaced by $(Q^2/B^{1-\al})^{n-1}$. Thus  (\ref{BQfrac}) is replaced by 
\[ Q^n(QP^{-2}(\log B)^2 + QB^{-\al(M-1)} + B^{(n-1)\al+1} Q^{-1})
\ll Q^{n+1},
\]
upon taking $\al=0$ and using $Q \gg B^{1/2}$. Ultimately, arguing in this way for $n=2$ leads to the choice $Q=B^{1/2} (\log B)^{1/2}$ and the outcome $S(F,B)\ll B^{n-1+1/2}(\log B)^{1/2}$, which is essentially no better than  (\ref{Cohen}), aside from the fact that we can remove the dependence on $\|F\|$ in the implicit constant.  In any case, Broberg's results (\ref{Bro_n2}) and (\ref{Bro_n3})  supercede the outcome of the methods of this paper for $n=2,3$.
\end{remark}

\section{Proof of Proposition \ref{prop_G_no_linear}}\label{sec_G}
In this section we prove the critical Proposition \ref{prop_G_no_linear} that allows us to deduce all factors in $G(0,\Ubf)$ have at least degree 2, so that we can apply the nontrivial bounds of Heath-Brown and Pila in (\ref{HBP_bounds}). We thank Per Salberger for suggesting the following strategy to prove the proposition.

 Let $n \geq 3$. Suppose to the contrary  that $G(0,\bfU)$ contains a linear factor, that is, 
\beq\label{G_has_linear}
G(0,\bfU)=L(\bfU)\tilde{H}(\bfU)
\eeq
for some linear form $L.$ Then by a linear change of variables we can reduce to the case in which we may assume that $L(\bfU)=U_1$, and conclude that 
\[G(0,\bfU)=U_1 H(\bfU)\]
for some homogeneous polynomial $H$.
Then \emph{any} point $(0,0,u_2,\ldots,u_n) \in \{U_Y=U_1=0\} \subset \mathbb{P}^n$  satisfies $G(0,\bfU)=0$ and thus defines a tangent hyperplane to $V(F(Z^e,\Xbf)) \subset \mathbb{P}^n$, given by 
\[ u_2X_2 + \ldots + u_nX_n=0.\]
In particular, for all $[u_2:\ldots:u_n] \in \mathbb{P}^{n-2}$, this hyperplane contains the line $\ell$ given by $X_2=\ldots =X_n =0$ in $\mathbb{P}^n$. We note that this line $\ell$ is not contained in $V(F(Z^e,\Xbf)),$ since for example 
in the coordinates $[U_Y:U_1:U_2:\ldots:U_n]$ we see that the point $[1:0:0:\ldots:0] \in \ell$ but $[1:0:0:\ldots:0] \not\in V$, since in the definition of $F$ the coefficient of $Z^{mde}$ is 1.
Thus under the assumption (\ref{G_has_linear}) we have shown that the generic 
hyperplane through $\ell$ is tangent to $V(F(Z^e,\Xbf))$. We will see this is impossible, and our assumption  (\ref{G_has_linear}) is false (so that Proposition \ref{prop_G_no_linear} is verified), by the following proposition.

\begin{prop}\label{prop_Salberger}
 Let $n \geq 3.$ Let $X \subset \mathbb{P}^n$ be a nonsingular hypersurface and let $\ell$ be a line not contained in $X$. Then the generic hyperplane in $\mathbb{P}^n$ containing $\ell$ is not tangent to $X$.
\end{prop}

Let $X$ be given as in the proposition. Without loss of generality we can make a change of coordinates so that 
\[\ell=\{X_2= \ldots =X_n=0\}.\]
Let $F \in \C[X_0,X_1,\ldots,X_n]$ be  such that $X=\{F=0\}$, and let $D$ denote the degree of $F$.
 Our strategy is to construct the blow-up of $X$ along the zero-dimensional subvariety $Z \subset X$, where we define
 \[ Z=\ell \intersect X \subset \mathbb{P}^n.\] 
 Under the hypothesis that $\ell$ is not contained in $X$, then $\deg Z \leq D.$ We also define the open set
 \[U :=X \setminus Z.\]
 To prove the proposition, we first notice that we can parametrize the hyperplanes containing $\ell$ in $\mathbb{P}^n$ by points in $\mathbb{P}^{n-2}$ using the map
 \[
 \begin{matrix}
 \mathbb{P}^{n-2} &\rightarrow &\{H\subset\mathbb{P}^{n}:\deg H=1,\text{ }\ell\subset H\}\\
 [v_{2}:\ldots:v_{n}] &\mapsto & \{v_{2}X_{2}+\ldots+v_{n}X_{n}=0\}.
 \end{matrix}
 \]
 Thus, it will suffice to show that there exists an open set $V \subset \mathbb{P}^{n-2}$ such that for all $\vbf=[v_2:\ldots: v_n] \in V,$ 
 \[ X \intersect \{v_2X_2 + \ldots + v_nX_n=0\}\]
 is smooth, so that in particular the hyperplane $\{v_2X_2+\ldots +v_nX_n=0\} \subset \mathbb{P}^n$ is not tangent to $X$. We will prove this in two steps, first focusing on the intersection of the hyperplane with the open set $U=X \setminus Z,$ and then focusing on the intersection of the hyperplane with the finite set of points in $Z$. In agreement with the citations we apply in what follows, from now on we will use the terminology ``regular'' for a scheme instead of ``smooth.'' For a nonsingular hypersurface such as $X$, these notions are identical  by the Jacobian criterion \cite[Ch. 4 Thm. 2.19, Ex. 2.10]{Liu02}; more generally, the notions are equivalent for any algebraic variety over a perfect field, and in particular over $\C$ \cite[Ch. 4 Cor 3.33]{Liu02}.

  Define a rational map $\varphi : X \dashrightarrow \mathbb{P}^{n-2}$ given by 
 \[\varphi: [X_0:X_1:X_2:\ldots:X_n] \mapsto [X_2: \ldots : X_n].\]
 This is a  regular map on $U$. We claim that there exists a projective variety $\tilde{Y}$ and two morphisms $\pi: \tilde{Y}\rightarrow X$, and $\tilde{\varphi}:\tilde{Y}\rightarrow \mathbb{P}^{n-2}$ such that
\begin{itemize}
\item[$i)$] The diagram
\[
  \begin{tikzcd}
    \tilde{Y}  \arrow{dr}{\tilde{\varphi}}  \arrow[swap]{d}{\pi} \\
   X  \arrow[dashrightarrow]{r}{\varphi} & \mathbb{P}^{n-2}
  \end{tikzcd}
\]
is commutative.
\item[$ii)$] the morphism $\pi$ restricts to an isomorphism $\pi:\pi^{-1}(U)\rightarrow U$.
\item[$iii)$] the projective variety $\tilde{Y}$ is regular.
\end{itemize} 
Let us assume this claim for now and see how to conclude the proof of the proposition. Since $\tilde{Y}$ is regular, we can apply Kleiman's Bertini theorem \cite[Ch. III Cor. 10.9]{Har77} to the morphism $\tilde{\varphi}:\tilde{Y}\rightarrow\mathbb{P}^{n-2}$, and deduce that given a generic hyperplane $H \subset \mathbb{P}^{n-2},$ $\widetilde{\varphi}^{-1}(H) \subseteq \widetilde{Y}$ is regular. Let us fix one of these generic hyperplanes, and call it 
\[ H = \{ u_2 X_2 + \ldots + u_n X_n=0\} \subset \mathbb{P}^{n-2}.\]
 By the choice of $H$, $\widetilde{\varphi}^{-1}(H) \intersect \pi^{-1}(U)$ is nonsingular.  Recall that  $\pi$ is an isomorphism when restricted to the open set $\pi^{-1}(U)$. Thus we also learn that
\begin{multline*}
\pi (\widetilde{\varphi}^{-1}(H) \intersect \pi^{-1}(U))  =
\pi (\widetilde{\varphi}^{-1}(H)  ) \intersect U 
= \varphi^{-1}(H) \intersect U \\= \{ [x_0:x_1:x_2:\ldots:x_n] \in U : u_2x_2 + \ldots + u_nx_n=0\}
\end{multline*} 
is regular.
Since such $H$ are generic in $\mathbb{P}^{n-2},$ we conclude that there is an open set $V_1 \subset \mathbb{P}^{n-2}$ such that for all $\bfv =[v_2:\ldots:v_n]\in V_1,$ the intersection 
\[ U \intersect \{ v_2X_2 + \ldots + v_nX_n=0\} \] is regular.

Let us next focus on the intersection of the hyperplane with the set $Z$. For any $P\in Z$, a hyperplane $\{v_{2}X_{2}+\ldots+v_{n}X_{n}=0\}$ with $[v_2: \ldots:v_n]\in \mathbb{P}^{n-2}$ is tangent to $X$ at $P$ if the Jacobian matrix at $P$,
\[
J_{\bfv}(P)=\begin{pmatrix}
\frac{\partial F}{\partial{X_{0}}}(P) & 0\\
\frac{\partial F}{\partial{X_{1}}}(P) & 0\\
\frac{\partial F}{\partial{X_{2}}}(P) & v_{2}\\
\vdots & \vdots\\
\frac{\partial F}{\partial{X_{n}}}(P) & v_{n}
\end{pmatrix},
\]
has rank $\leq 1$. From this it is clear that if either $\frac{\partial F}{\partial{X_{0}}}(P)\neq 0$ or $\frac{\partial F}{\partial{X_{1}}}(P)\neq 0$ then $\rank J_{\bfv}(P)=2$ for any $\bfv\in\mathbb{P}^{n-2}$. On the other hand, if $\frac{\partial F}{\partial{X_{0}}}(P)=\frac{\partial F}{\partial{X_{1}}}(P)=0$ then $\rank _{\bfv}(P)\leq 1$ if and only if $\bfv=[\frac{\partial F}{\partial{X_{2}}}(P):\ldots:\frac{\partial F}{\partial{X_{n}}}(P)]$ since we are assuming that $X$ is a nonsingular hypersurface. For each $P\in Z$ we define
\[
C_{P}=
\begin{cases}
\{[\frac{\partial F}{\partial{X_{2}}}(P):\ldots:\frac{\partial F}{\partial{X_{n}}}(P)]\}           &\text{if $\frac{\partial F}{\partial{X_{0}}}(P)=\frac{\partial F}{\partial{X_{1}}}(P)=0$,}\\
\emptyset                                                                                       &\text{otherwise}.
\end{cases}
\]
If we define $V_{P}=\mathbb{P}^{n-2}\setminus C_{P}$, it follows that for any $\bfv\in V_{P}$ the intersection
\[
X \intersect \{ v_2X_2 + \ldots + v_nX_n=0\}
\]
is regular at $P$. 

Finally consider the set 
\[V=V_{1}\cap\bigcap_{P\in Z}V_{P}.\]
Since $\deg Z\leq D$, then $V$ is a non-empty open subset of $\mathbb{P}^{n-2}$. For each $\vbf \in V$, the hyperplane  $v_2x_2 + \ldots + v_nx_n=0$ contains $\ell$, and \[\{v_2X_2+\ldots+v_nX_n=0\}\intersect (U\union Z)=\{v_2X_2+\ldots+v_nX_n=0\}\intersect X\]
is regular, or equivalently, nonsingular; thus $\{v_2X_2+\ldots+v_nX_n=0\}$ is not tangent to $X.$ This completes the proof of Proposition \ref{prop_Salberger}, except for the proof of properties (i), (ii), and (iii) in the claim. 

We now prove the claim of properties (i), (ii) and (iii). From the rational map $\varphi : X \dashrightarrow \mathbb{P}^{n-2}$ given by 
 \[\varphi: [X_0:X_1:X_2:\ldots:X_n] \mapsto [X_2: \ldots : X_n],\]
we consider the graph $\Ga=\Ga_\varphi$ of the map $\varphi$, 
\[ \Ga = \{ (\x, \varphi(\x)) : \x \in U\} \subset X \times \mathbb{P}^{n-2}.\]
Define the Zariski closure $\widetilde{X} = \overline{\Ga} \subset X \times \mathbb{P}^{n-2}.$
 Define the projection map $\pi' : \widetilde{X} \maps X$ acting by $(\x,\varphi(\x)) \maps (\x).$
Then the blow-up is $\tilde{X}$ along with a morphism $\varphi'$ such that 
\[
  \begin{tikzcd}
    \tilde{X}  \arrow{dr}{\varphi'}  \arrow[swap]{d}{\pi'} \\
   X  \arrow[dashrightarrow]{r}{\varphi} & \mathbb{P}^{n-2}
  \end{tikzcd}
\] 
is a commutative diagram (see e.g. \cite[Ch. 7 p. 82]{Har92}). Moreover, from the definition of the blow-up it follows that   $\pi'$ restricts to an isomorphism $\pi':(\pi')^{-1}(U)\rightarrow U$, i.e. $\tilde{X}$ satisfies properties (i) and (ii),  but it might be singular. To resolve this, we apply Hironaka's resolution of singularities: as a consequence of \cite[Theorem $1$]{Hir64a} (see also \cite[P. 112]{Hir64a}), there is a projective variety $\tilde{Y}$ and a morphism $f:\tilde{Y}\rightarrow \tilde{X}$ such that $f$ is an isomorphism when restricted to the inverse image $f^{-1}(V)$ of the open set $V$ of the regular points of $\tilde{X}$, and such that $\tilde{Y}$ is regular. Then the claim follows by taking $\pi=\pi'\circ f$, $\widetilde{\varphi}=\varphi'\circ f$ and observing that $(\pi')^{-1}(U)\subset V$.

\section{Concluding arguments}\label{sec_concluding}

In \S \ref{sec_badbad} we proved that the contribution of the bad-bad terms to the sieve is $\ll Q^{n}.$ We now turn to analyzing the contributions of the other types, as defined in Definition \ref{dfn_types}. We will treat these in three sections; in each case we apply the relevant bound for $|g(\ubf,pq)|$ from Proposition \ref{prop_g_sum} and the bound (\ref{WB_bound_sum}) for $\hat{W}$. Once we have treated these cases, we proceed in \S \ref{sec_final} to choose the parameter $Q$, and conclude the proof of Theorem \ref{thm_main}.

\subsection{Zero-type cases}\label{sec_zero}
We first consider any case in which $\ubf$ is zero-type modulo $p$, divided into cases according to whether $\ubf$ is zero-type, good, or bad modulo $q$. The contribution of the first case (upon setting $\ubf = pq\vbf$ and applying (\ref{WB_bound_sum})) is 
\[
\frac{1}{P^2 Q^{2n}}\sum_{\substack{p,q\in\mathcal{P}\\p\neq q}}\sum_{\substack{\bfu\in\mathbb{Z}^{n}\\\bfu\text{ zero mod }p\\\bfu\text{ zero mod }q}}\left|\hat{W}\left(\frac{\bfu}{pq}\right) g(\bfu,pq)\right|
\ll 
\frac{Q^{2n-1}}{P^{2}Q^{2n}}\sum_{\substack{p,q\in\mathcal{P}\\p\neq q}}\sum_{ \bfv\in\mathbb{Z}^{n}}\left|\hat{W} (\bfv)\right| \ll B^nQ^{-1}.
\]
The contribution of the second case (upon setting $\ubf = p\vbf$, applying (\ref{WB_bound_sum}) with $L=Q<B$) is
\[
\frac{1}{P^2 Q^{2n}}\sum_{\substack{p,q\in\mathcal{P}\\p\neq q}}\sum_{\substack{\bfu\in\mathbb{Z}^{n}\\\bfu\text{ zero mod }p\\\bfu\text{ good mod }q}}\left|\hat{W}\left(\frac{\bfu}{pq}\right) g(\bfu,pq)\right|
\ll 
\frac{Q^{n-1/2}Q^{n/2}P^2}{P^2Q^{2n}} \sum_{\vbf \in \Z^n}\left|\hat{W}\left( \frac{\vbf}{Q}\right)\right|
\ll B^nQ^{-n/2-1/2}.
\]
The contribution of the third case (upon  setting $\ubf = p\vbf$, applying (\ref{WB_bound_sum}) with $L=Q<B$) is
\[\frac{1}{P^2 Q^{2n}}\sum_{\substack{p,q\in\mathcal{P}\\p\neq q}}\sum_{\substack{\bfu\in\mathbb{Z}^{n}\\\bfu\text{ zero mod }p\\\bfu\text{ bad mod }q}}\left|\hat{W}\left(\frac{\bfu}{pq}\right) g(\bfu,pq)\right|
\ll 
\frac{Q^{n-1/2}Q^{n/2+1/2}P^2}{P^2Q^{2n}} \sum_{\vbf \in \Z^n}\left|\hat{W}\left( \frac{\vbf}{Q}\right)\right|
\ll B^nQ^{-n/2}.
\]
As long as $n \geq 2$, all these cases contribute at most $\ll B^nQ^{-1}$ to the sieve, which is acceptable.

\subsection{Good-good case}
The contribution to the sieve from the good-good case is:
\[
\frac{1}{P^2 Q^{2n}}\sum_{\substack{p,q\in\mathcal{P}\\p\neq q}}\sum_{\substack{\bfu\in\mathbb{Z}^{n}\\\bfu\text{ good mod }p\\\bfu\text{ good mod }q}}\left|\hat{W}\left(\frac{\bfu}{pq}\right) g(\bfu,pq)\right|
\ll \frac{Q^nP^2}{P^2 Q^{2n}} \sum_{\bfu\in\mathbb{Z}^{n}}\left|\hat{W}\left(\frac{\bfu}{Q^2}\right)  \right|
    \ll Q^n,
\]
after applying (\ref{WB_bound_sum}) with $L=Q^2>B$, since under the assumption (\ref{kappa_assp}), $\kappa \geq 1/2.$
\subsection{Good-bad case}\label{sec_good_bad}
The contribution to the sieve from the good-bad case is
\beq\label{good-bad}
\frac{1}{P^2 Q^{2n}}\sum_{\substack{p,q\in\mathcal{P}\\p\neq q}}\sum_{\substack{\bfu\in\mathbb{Z}^{n}\\\bfu\text{ good mod }p\\\bfu\text{ bad mod }q}}\left|\hat{W}\left(\frac{\bfu}{pq}\right) g(\bfu,pq)\right|
\ll 
\frac{Q^{n+1/2}}{P^2 Q^{2n}}\sum_{p \in \Pcal}\sum_{q \neq p \in \Pcal}\sum_{\substack{\bfu\in\mathbb{Z}^{n}\\ \bfu\text{ bad mod }q}}\left|\hat{W}\left(\frac{\bfu}{pq}\right)  \right|.
\eeq
Here we proceed by imitating the key step from \S \ref{sec_badbad} for the bad-bad case, and sum over $q$ before summing over $\ubf$. We again define $G(U_Y,\Ubf)$ as in (\ref{G_dfn}), and let $R(\ubf)$ denote the resultant (\ref{R_dfn}), so that 
\[
\sum_{p \in \Pcal} \sum_{\substack{\bfu\in\mathbb{Z}^{n}\\ G(0,\ubf)\neq 0}}\sum_{\substack{q \neq p \in \Pcal\\ \bfu\text{ bad mod }q}}\left|\hat{W}\left(\frac{\bfu}{pq}\right)  \right|
    \ll P \sum_{\substack{\bfu\in\mathbb{Z}^{n}\\ G(0,\ubf)\neq 0}}\left|\hat{W}\left(\frac{\bfu}{Q^2}\right)  \right| \om(R(\ubf))
   \ll_{n,m,e,d} P  Q^{2n} \log B ,
\]
with an implied constant independent of $\|F\|$ (in the first case of Lemma \ref{lemma_F_small}),    by arguing as in the proof of (\ref{badbad_unonzero}).  
    
 Notice that in the good-bad case, we do not need to consider a possible contribution from those $\ubf$ for which  $G(0,\ubf)=0$: when $G(0,\ubf)=0$, then all $q$ have the property that $\ubf$ is bad for $q$, whereas by definition in the good-bad case, $\ubf$ is good for at least one prime.
In total, the contribution to the sieve from the good-bad case is thus
\[
     \frac{Q^{n+1/2}}{P^2Q^{2n}} \cdot P Q^{2n}(\log B) 
     \\
 \ll  Q^{n+1/2}P^{-1}(\log B) \ll Q^n,
\]
since $Q=B^\kappa$ for some $1/2 \leq \kappa \leq 1$ and under our acting assumption (\ref{Q_lower}), by (\ref{P_lower}), $P \gg Q/\log Q.$
Thus we can conclude that the total contribution of the good-bad case (\ref{good-bad}) of the sieve is $\ll  Q^{n}$, with an implied constant independent of $\|F\|$ (in the first case of Lemma \ref{lemma_F_small}).

\subsection{Final conclusion of the sieve, and choice of parameters}\label{sec_final}

We now assemble all the terms of the main sieve term in (\ref{sieve_Tpq}): we can conclude that 
\beq\label{T_final}
\frac{1}{P^{2}}\sum_{\substack{p,q\in\mathcal{P}\\p\neq q}}|T(p,q)|\ll B^nQ^{-1} + Q^{n}.
\eeq
The first term is from all zero-type cases, and the last term includes the good-good, good-bad, and bad-bad cases.
We apply this in the sieve lemma, along with the bound (\ref{S_simple}) for the two simple terms in the sieve, to conclude that (in the first case of Lemma \ref{lemma_F_small}) our counting function admits the bound
\beq\label{S_final_bound}
\mathcal{S}(F,B)\ll_{n,m,e,d}  \left( B^{n-1} + B^nP^{-1}+B^nQ^{-1} + Q^{n} \right) \ll   \left( B^nP^{-1} + Q^n \right).
\eeq
 Choose 
\beq\label{Q_choice}
Q=B^{n/(n+1)}(\log B)^{1/(n+1)}.
\eeq
The requirement (\ref{kappa_assp}) is met for all $n \geq 3$. (If $n=2$, then this argument leads to the choice $Q \approx B^{2/3}$, which does not suffice to prove sufficient decay in the bad-bad case; see Remark \ref{remark_badbad_n2}.)
Recall from (\ref{Q_lower}) and (\ref{P_lower}) that 
\[P = |\Pcal| \gg_{m,e,d} Q(\log Q)^{-1} \gg_{n,m,e,d} B^{\frac{n}{n+1}}(\log B)^{-\frac{n}{n+1}}
\]
as long as 
\beq\label{Q_require}
Q \gg_{m,e,d}(\log \|F\|)(\log \log \|F\|).
\eeq
Recall also that we require
$P \gg_{m,e,d} \max\{ \log \|f_d\|, \log B\}$ in Lemma \ref{lemma_sieve}. Certainly the first condition is satisfied under the assumption (\ref{Q_require}). The second condition is satisfied for $Q$ as in (\ref{Q_choice}) for all $B\gg_n 1.$

To meet the requirement (\ref{Q_require}) for $Q$ as chosen in (\ref{Q_choice}), it suffices to require that
\[ B \gg_{m,e,d} (\log \|F\| \log \log \|F\|)^{\frac{n+1}{n}}.\]
For such $B$,  the conclusion of the sieve process in (\ref{S_final_bound}) shows that  
\[
\mathcal{S}(F,B)\ll_{n,m,e,d}   B^{n-1+\frac{1}{n+1}} (\log B)^{\frac{n}{n+1}},
\]
where the implicit constant is independent of $\|F\|.$
This suffices for Theorem \ref{thm_main}.
Finally, for all $B \ll_{m,e,d}(\log \|F\| \log \log \|F\|)^{\frac{n+1}{n}}$, we apply the trivial bound 
\begin{multline*}
\Scal(F,B) \ll_n B^n \ll_{n,m,e,d} (\log \|F\| \log \log \|F\|)^{n+1} \ll (\log \|F\|)^{n+2} \\
\ll_{n,m,d,e} (\log B)^{n+2} \ll_n B^{n-1+\frac{1}{n+1}} (\log B)^{\frac{n}{n+1}} .
\end{multline*}
Here we applied the fact from Lemma \ref{lemma_F_small} that in the case it remains to prove Theorem \ref{thm_main}, $\|F\|\ll B^{(mde)^{n+2}}$ so that $ \log \|F\| \ll_{n,m,d,e} \log B$. This completes the proof of Theorem \ref{thm_main}.

\section*{Funding}

The first author has been supported by FWF grant P 32428-N35. The second author   has been partially supported by NSF DMS-2200470 and NSF
CAREER grant DMS-1652173, a Sloan Research Fellowship, and a Joan and Joseph Birman Fellowship. The authors thank the Hausdorff Center for Mathematics for hosting a productive collaboration visit and the RTG DMS-2231514; the second author thanks HCM for hosting visits as a Bonn Research Chair.

\bibliographystyle{alpha}

\bibliography{NoThBibliography}


 \pagebreak
 \pagestyle{plain}
   \setcounter{page}{1}

 \setcounter{section}{0}

  \setcounter{equation}{0}
  \renewcommand{\theequation}{\arabic{equation}}

 \begin{center}
 {\bf Correction to ``Application of a polynomial sieve: beyond separation of variables''}\\
 \vspace{.5cm}
 Correction as published in \emph{Algebra \& Number Theory} (2026)\\
  \vspace{.5cm}
Dante Bonolis and  Lillian B. Pierce 
  \end{center}
  
 \vspace{1cm}

Fix an integer $m \geq 2$ and integers $d, e \geq 1$. Consider a polynomial
\[
F(Y,\bX)=Y^{md}+Y^{m(d-1)}f_{1}(\bX)+\cdots+Y^mf_{d-1}(\bX)+f_{d}(\bX),
\]
in which for each $1 \leq i \leq d$,  $f_i \in \Z[X_1,\ldots, X_n]$ is a form with $\deg f_{i}=m\cdot e\cdot i$ (and $f_d \not\con 0$).
Define
\[
N(F, B):=\#\{\bx\in [-B,B]^{n} \cap \Z^n: \exists y\in\mathbb{Z}\text{ such that } F(y,\bx)=0\}.
\]
Fix $n\geq 3$, and suppose the weighted hypersurface $V(F(Y,\bX)) \subset \mathbb{P}(e,1,\ldots,1)$ defined by $F(Y,\bX)=0$ is nonsingular over $\C$.  Let $\|F\|$ denote the maximum absolute value of any coefficient of the polynomial $F$; it is no loss of generality below to assume that $\|F\| \geq 3$ and $B \geq 3$.
Theorem 1.1 of \cite{BonPie24} states that  under the above hypotheses,
\beq\label{uniform}
N(F,B) \ll_{n,m,d,e} B^{n-1+\frac{1}{n+1}} (\log B)^{\frac{n}{n+1}}
\eeq
with an implied constant that can depend on $n,m,d,e$, but is independent of $\|F\|$.
Here we correct this to the statement:

{\bf Theorem 1.1'}:  Under the above hypotheses, for some positive integer $h(n)$,
\beq\label{polylog}
N(F,B) \ll_{n,m,d,e} (\log \|F\|)^{h(n)} B^{n-1+\frac{1}{n+1}}(\log B)^{\frac{n}{n+1}}.
\eeq

The bound stated in Theorem 1.1' is the direct outcome of the polynomial sieve, which is correctly proved in the main argument of the original paper \cite{BonPie24}; we briefly demonstrate in \S \ref{sec_tracking_dependence}   how to track the dependence on $\|F\|$. 

The original paper claims that (\ref{uniform}) can be obtained because (\ref{polylog}) can be upgraded to (\ref{uniform}) by an application of  Lemma 2.1 in \cite{BonPie24}.
  But the proof of Lemma 2.1  contains a gap, so the lemma is not valid and it cannot be applied. 
  
 Lemma 2.1   considers a hypersurface $V(G(Y,\bX))\subset\mathbb{P}(e,1,\ldots,1)$, defined by  
\[G(Y,\bX)=Y^D + Y^{D-1}f_1(\bX) + \cdots + Yf_{D-1}(\bX)+ f_D(\bX)\]
with each $f_i \in \Z[X_1,\ldots,X_n]$ a form of $\deg f_i = i \cdot e$, for fixed $D,e \geq 1$ and $n \geq 1$.  Lemma 2.1 claims that if $f_D \not\con 0$ and the weighted hypersurface $V(G(Y,\bX))\subset\mathbb{P}(e,1,\ldots,1)$ is absolutely irreducible, then either
\beq\label{NG_false}
N(G, B):=\#\{\bx\in [-B,B]^{n} \cap \Z^n: \exists y\in\mathbb{Z}\text{ such that } G(y,\bx)=0\}\ll_{n,D,e} B^{n-1}
\eeq
or $
\|G\|\ll B^{(De)^{n+2}}.
$ 
Here we correct this to the statement:

{\bf Lemma 2.1'}:  Under the above hypotheses, either 
\[
N'(G, B):=\#\{\bx\in [-B,B]^{n} \cap \Z^n: \exists y\in [-B^e,B^e]\cap \mathbb{Z}\text{ such that } G(y,\bx)=0\}\ll_{n,D,e} B^{n-1}
\]
or $
\|G\|\ll B^{(De)^{n+2}}.
$
The conclusion of Lemma 2.1', when applied to the polynomial $F(Y,\bX)$, is not useful to upgrade (\ref{polylog}) to (\ref{uniform}), since it refers to a modified counting function.
 The essential distinction is that $N'(G,B)$ additionally restricts $y$ to the interval $[-B^e,B^e]$ independent of $\|G\|$, whereas $N(G,B)$ does not. (Naively, for a given $\bx$ lying in the set counted by $N(G,B)$, if $y$ solves $G(y,\bx)=0$,   $|y|$ could be as large as $\|G\|^{1/D}B^e$.) 
In   \S \ref{sec_pinpoint} of this correction, we explicitly describe the gap in the proof of Lemma 2.1, and also indicate how to prove Lemma 2.1'.

\section{Proof of Theorem 1.1': tracking dependence on $\|F\|$}\label{sec_tracking_dependence}
For clarity, we verify here that the dependence in (\ref{polylog}) is only polylogarithmic in $\|F\|$, as a consequence of the argument already presented in the main body of \cite{BonPie24}; all equation numbers and section numbers refer to that reference. To do so, we now indicate all the places in \cite{BonPie24} with dependence on $\|F\|$. First, the sieving set must consist of primes sufficiently large with respect to $\|F\|$, as seen in two instances.  Equation (1.22) of Lemma 1.2 requires that $P=|\mathcal{P}| \gg_{m,e,d} \log \|F\|$. Equation (4.4) of  \S 4.4 requires $Q \gg_{m,d,e} (\log \|F\|)( \log \log \|F\|)$; this ensures that the previous  condition holds. Second, dependence on $\|F\|$ enters the argument of the polynomial sieve in order to control for how many primes $p$ a vector $\bfu$ can be ``locally bad'' (that is $G(0,\bfu) \neq 0$ as an integer but $p | G(0,\bfu)$).   Equations (5.7) and (5.10) show  
\[|\{p :\text{$\mathbf{u}$ is bad modulo $p$}\}| \leq \om (R(\mathbf{u}))\ll_{n,m,e,d} \log (\|F\| \|\mathbf{u}\|),\]
and the factor of $\log \|F\|$ appearing here will affect the bound proved for Equation (5.12); it will not be relevant for Equation (5.13). 
In \S 5.2.1 to bound Equation (5.12), we apply $\om(R(\bu))^2\ll (\log (\|F\|B))^2$, replacing the statement  $\om(R(\bu))^2\ll(\log B)^2$ as applied in the paper. Consequently, Equation (5.12) now has right-hand side $\ll_{n,m,e,d} (\log \|F\|)^2 Q^{2n} (\log B)^2$. Carrying the factor $(\log \|F\|)^2$ through the analysis of the bad-bad contribution in \S 5.2.3 finally shows Equation (5.22) now with right-most side $\ll(\log \|F\|)^2Q^n$. In \S 7.3, the good-bad contribution also carries one factor of $\om(R(\bu))\ll\log( \|F\| B)$, so the good-bad contribution is bounded by $\ll (\log \|F\|)Q^n$ in total.   Thus the final outcome of the polynomial sieve, Equation (7.2), holds with right-hand side $\ll B^n Q^{-1} + (\log \|F\|)^2Q^n$. 
Arguing exactly as in \S 7.4 then shows that for all 
\[ B \gg_{m,e,d} (\log \|F\| \log \log \|F\|)^{\frac{n+1}{n}},\]
the choice $Q=B^{n/(n+1)}(\log B)^{1/(n+1)}$ satisfies the requirement $Q \gg_{m,e,d}(\log \|F\|)(\log \log \|F\|)$, and  the conclusion of the sieve process is that 
\[
N(F,B) \ll_{n,m,e,d}  (\log \|F\|)^2 B^{n-1+\frac{1}{n+1}} (\log B)^{\frac{n}{n+1}}.
\]
Finally, for all $B \ll_{m,e,d}(\log \|F\| \log \log \|F\|)^{\frac{n+1}{n}}$,  apply the trivial bound 
\begin{multline*}
N(F,B) \ll_n B^n \ll_{n,m,e,d} (\log \|F\| \log \log \|F\|)^{n+1} \ll (\log \|F\|)^{n+2} \\
\ll_{n,m,d,e} (\log \|F\|)^{n+2}   B^{n-1+\frac{1}{n+1}} (\log B)^{\frac{n}{n+1}} .
\end{multline*}
This verifies (\ref{polylog}).

  \section{Proof of Lemma 2.1', and the gap in the proof of Lemma 2.1}\label{sec_pinpoint}

 We now pinpoint the gap in the proof presented in \cite[Lemma 2.1]{BonPie24} to control $N(G,B)$ as defined in (\ref{NG_false}), and specify how the proof successfully controls $N'(G,B)$ as defined in Lemma 2.1' of this correction. Recall the matrix $\mathbf{C}$ used in the proof method for Lemma 2.1,  constructed by 
\[
\mathbf{C}=(\bfv_{i}^{\bfe})_{\substack{1\leq i\leq N\\ \bfe\in\mathcal{E}}}.
\] 
Here $B\geq 1$ is fixed, and  $\{\bfv_{1},\ldots, \bfv_{N}\}$ enumerate the   solutions to $G(Y,\bX)=0$, in coordinates $(y,x_{1},\ldots,x_{n})$, with  each of $x_1,\ldots,x_n$ lying in $[-B,B] \cap \Z$ and \emph{no imposed constraint on the size of $y \in \Z$}.  
The proof correctly constructed  a nonzero vector $\bfb \in \Z^{|\mathcal{E}|}$ in the nullspace of $\mathbf{C}$  with entries that are $(|\mathcal{E}|-1) \times (|\mathcal{E}|-1)$ minors of $\mathbf{C}$, and with the property that $\|G\| \leq |\bfb|$. In particular, if we let $C_{\max}$ represent the maximum absolute value of any entry in $\mathbf{C}$ then it is true that $\|G\| \leq |\bfb| \ll_{|\mathcal{E}|} C_{\max}^{|\mathcal{E}|}$. The proof of Lemma 2.1 effectively claimed that $C_{\max} \ll B^{De}$, independent of $\|G\|$, from which it would follow
 $|\bfb|\ll B^{De|\mathcal{E}|}\ll B^{(De)^{n+2}}$. But the claim $C_{\max} \ll B^{De}$  is false. An entry   in $\mathbf{C}$ can for example be as big as $|y^D|  \approx \|G\|B^{De}$, which depends on $\|G\|$.  Thus   $C_{\max}$ cannot be bounded independent of $\|G\|$ a priori, and the strategy described to prove the lemma cannot guarantee the second outcome claimed in the lemma, namely $\|G\|\ll B^{(De)^{n+2}}$.

 However, the strategy described in Lemma 2.1 of \cite{BonPie24} succeeds to prove a dichotomy for the modified counting function $N'(G,B)$. For the modified counting function, the proof can proceed by constructing instead a matrix $ \mathbf{C}'=(\tilde{\bfv}_{i}^{\bfe})$, in which 
  $\{\tilde{\bfv}_{1},\ldots, \tilde{\bfv}_{N'}\}$ enumerate the solutions $(y,x_1,\ldots,x_n)$   to $G(Y,\bX)=0$,   with  each of $x_1,\ldots,x_n$ lying in $[-B,B] \cap \Z$ \emph{and with  the additional constraint} $y \in [-B^e,B^e] \cap \Z$. 
In this setting, the construction outlined in \cite{BonPie24} correctly constructs  a nonzero vector $\bfb' \in \Z^{|\mathcal{E}|}$ in the nullspace of $\mathbf{C}'$  with entries that are $(|\mathcal{E}|-1) \times (|\mathcal{E}|-1)$ minors of $\mathbf{C}'$, and with the property that $\|G\| \leq |\bfb'|$. Now, if we let $C_{\max}'$ represent the maximum absolute value of any entry in $\mathbf{C}'$ then it is true that $C_{\max}' \ll B^{De}$, independent of $\|G\|$, and this leads to the conclusion $\|G\| \ll B^{(De)^{n+2}}$.
To summarize, the strategy of proof given for Lemma 2.1 in \cite{BonPie24} is   valid in settings in which all the variables under consideration are constrained by a box that depends only on $B$ and not $\|G\|$. 

The authors thank Katharine Woo for discussions on these topics.

\end{document}